\documentclass[11pt]{article}
\RequirePackage[aos,amsthm,amsmath,natbib]{imsartmine}

\usepackage{epsfig}
\usepackage{amsmath,amsfonts,bm,rotating,amssymb,amsthm,latexsym,subfigure}
\usepackage{longtable}
\usepackage{multirow}

\startlocaldefs

\renewcommand{\thefigure}{\thesection.\arabic{figure}}

\def\definedas{\stackrel{\Delta}{=}}

\def\:{:\,}

\newcommand{\bc}{\begin{center}}
\newcommand{\ec}{\end{center}}
\newcommand{\beq}{\begin{eqnarray}}
\newcommand{\eeq}{\end{eqnarray}}
\newcommand{\beqq}{\begin{eqnarray*}}
\newcommand{\eeqq}{\end{eqnarray*}}

\def\LK{Lipschitz-Killing}

\def\LKCs{Lipschitz-Killing curvatures\ }
\newcommand{\lips}{{\cal L}}
\def\EC{Euler characteristic\ }

\def\E{E}
\def\P{P}
\def\RN{{R}^N}

\def\real{{\cal R}}
\def\definedas{\stackrel{\Delta}{=}}
\def\:{:\,}
\def\p{\varphi}

\endlocaldefs

\begin{document}
\begin{frontmatter}

\title{Estimating Thresholding Levels for Random Fields via Euler Characteristics}
\runtitle{Thresholds via Euler characteristics}

\begin{aug}
\author{\fnms{Robert J.} \snm{Adler}\thanksref{t1,t3}
\ead[label=e1]{robert@ee.technion.ac.il}
\ead[label=u1,url]{robert.net.technion.ac.il}}
\author{\fnms{Kevin} \snm{Bartz}\thanksref{t1}
\ead[label=e2]{kevin@kevinbartz.com}
\ead[label=u2,url]{????}
}
\author{\fnms{S.\ C.} \snm{Kou}\thanksref{t1,t4}
\ead[label=e3]{kou@stat.harvard.edu}
\ead[label=u3,url]{www.people.fas.harvard.edu/$\sim$skou}}
\and 
\author{\fnms{Anthea} \snm{Monod}\thanksref{t2}
\ead[label=e4]{am4691@cumc.columbia.edu}
\ead[label=u4,url]{systemsbiology.columbia.edu/people/anthea-monod}}

\thankstext{t1}{Research supported in part by US-Israel Binational
Science Foundation, 2008262.}
\thankstext{t2}{Research supported in part by FP7-ICT-318493-STREP.}
\thankstext{t3}{Research supported in part by  ERC  2012 Advanced Grant 20120216.}
\thankstext{t4}{Research supported in part by NIH/NIGMS,R01GM090202 and NSF, DMS-1510446.}


\affiliation{Technion, Renaissance Technologies, Harvard and Columbia.}
\runauthor{Adler, Bartz, Kou and Monod}
\end{aug}



\begin{abstract}
We introduce Lipschitz-Killing curvature (LKC) regression, a new method to produce $(1-\alpha)$ thresholds for signal detection in random fields that does not require knowledge of the spatial correlation structure.  The idea is to fit  observed empirical Euler characteristics to the Gaussian kinematic formula via generalized least squares, which quickly and easily provides statistical estimates of the LKCs --- complex topological quantities that can  be extremely challenging to compute, both theoretically and numerically.  With these estimates, we can then make use of a powerful parametric approximation via Euler characteristics for Gaussian random fields to generate accurate $(1-\alpha)$ thresholds and $p$-values.  The main features of our proposed LKC regression method are easy implementation, conceptual simplicity, and facilitated diagnostics, which we demonstrate in a variety of simulations and applications.

\end{abstract}
\end{frontmatter}
\section{Introduction}

\label{SectionIntroduction}

\subsection{Pre-History}
The paper that you are (maybe) about to read was started in 2008, when RJA spend a semester at Harvard and talked a lot with  SCK.
During a lecture, RJA made the claim that estimating \LKCs (defined below) was a hard problem, at which stage SCK commented ``It doesn't seem all that hard.  Why can't you just treat it as a regression problem?"  At the time KB was a graduate student, listening to the same lecture, and all three got together to write the first version of this paper. Since RJA returned to the Technion, and KB graduated and moved to Renaissance Technologies, it took until 2011 for a final first version to be submitted for publication. A few months later it was rejected (from an excellent journal) but, since the authors still liked their paper they set about  preparing a revision, this time together with AM, who was then a postdoctoral fellow at the Technion. In 2014 the new version was submitted to a different (and even better) journal, only to be rejected since it was too methodological and not applied enough, but with some nice ideas for rewriting as long as the authors were prepared to rewrite everything in such a way that it made the referees happy and them miserable.

Now that AM has left the Technion and is doing other things, RJA is (happily) approaching retirement,  KB is well entrenched in industry,  and AM is also busy with new problems, it seems unlikely that this paper will ever get revised again. But, the four authors still like it, and it keeps getting cited (seemingly implying that others share their opinion) and so the time has come to put what there is  up on arXiv. Which is why what you are reading is a 2017 arXiv deposit of what is basically a 2011 paper.

\subsection{Motivation} 
Random field models are widely used in many scientific applications, including the description of spatial structures in environmental and epidemiological studies, the statistical analysis of brain images, and the modeling of the cosmic microwave background radiation, as well as other cosmological phenomena.  An important problem of interest common to most of these applications is the determination of threshold levels for the random field, which indicate that regions with values above the level are, in some sense, significant, while regions with values below are not.  Accurate determination of threshold levels for random fields faces a major challenge that the simpler setting of independent observations does not: the values of random fields are correlated in space.

\begin{figure}[!ht]
  \begin{tabular}{ccc}
    \textbf{Subject 1} & \textbf{Subject 1} & \textbf{Subject 1} \\ 
    \textbf{Same-speaker} & \textbf{Difference} & \textbf{Studentized Residual}
    \\
    \includegraphics[width=1.4in]{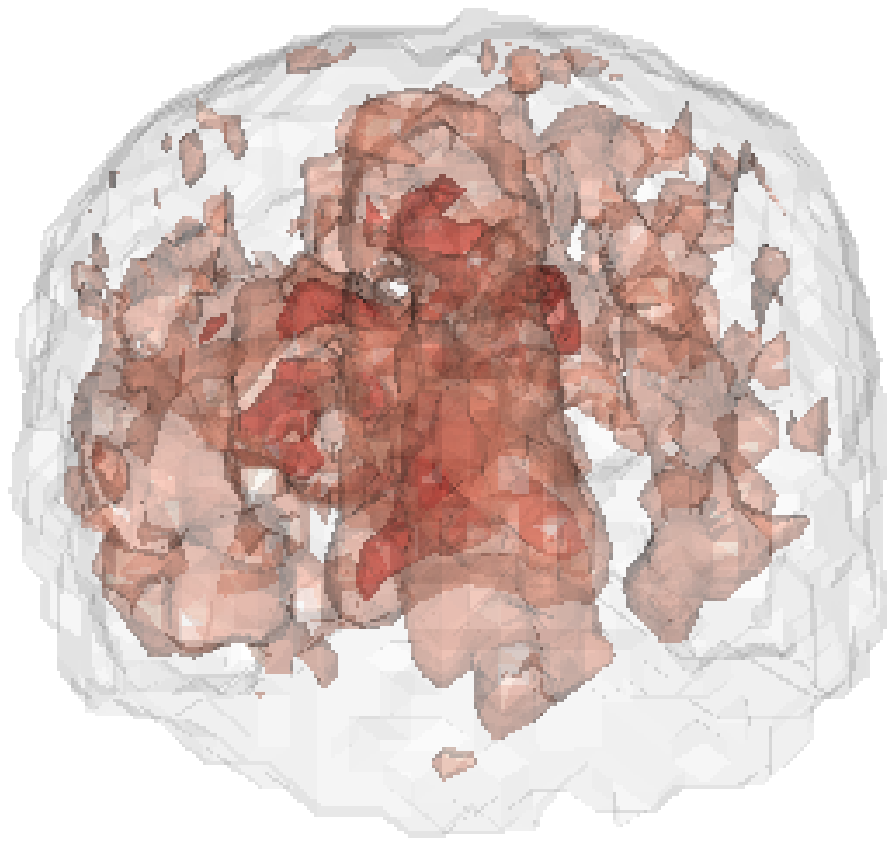} & %
    \includegraphics[width=1.4in]{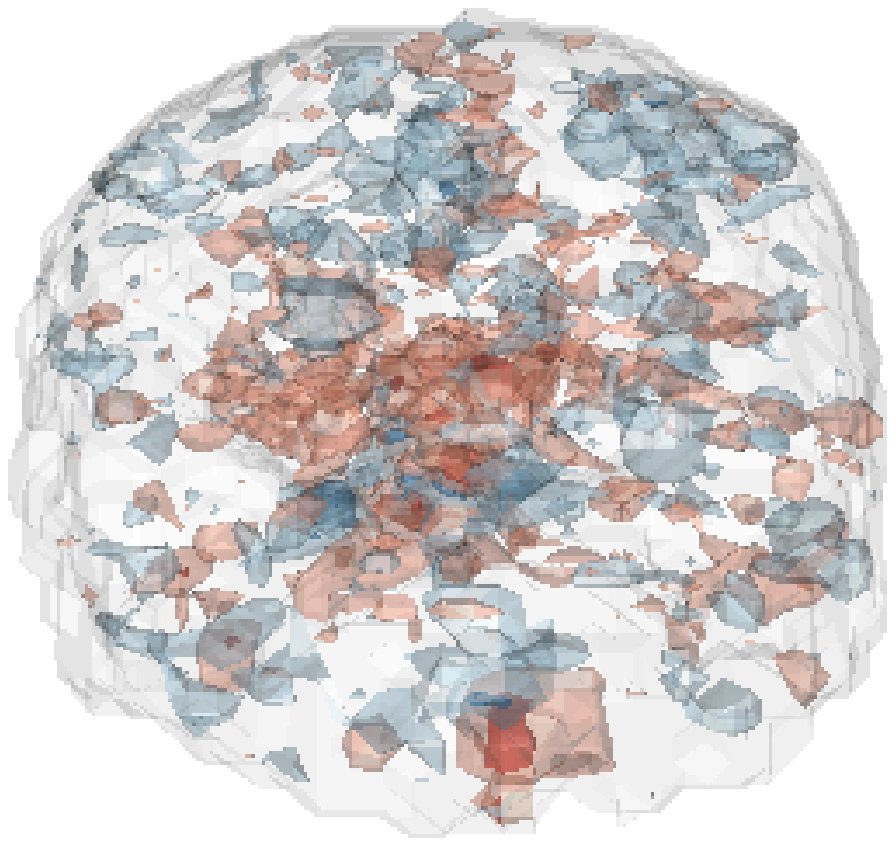} & %
    \includegraphics[width=1.4in]{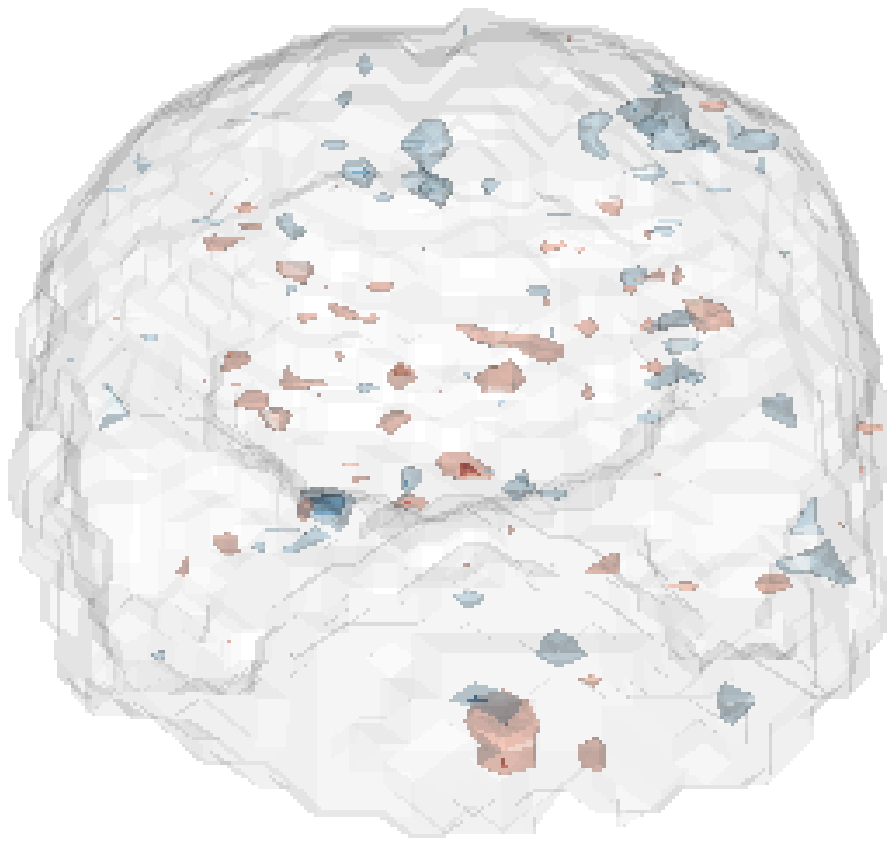} \\ 
    &  &  \\ 
    \textbf{Subject 1} & \textbf{All Subjects} & \textbf{All Subjects} \\ 
    \textbf{Different-speaker} & \textbf{Average Difference} & \textbf{$t$
      Statistic} \\ 
    \includegraphics[width=1.4in]{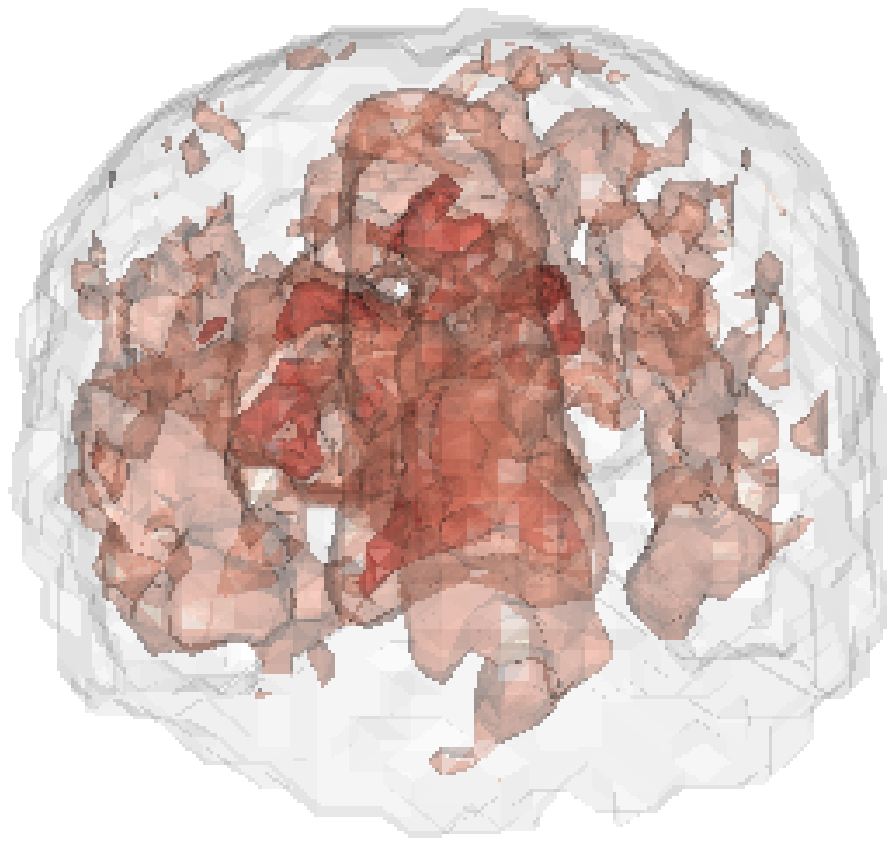} &
    \includegraphics[width=1.4in]{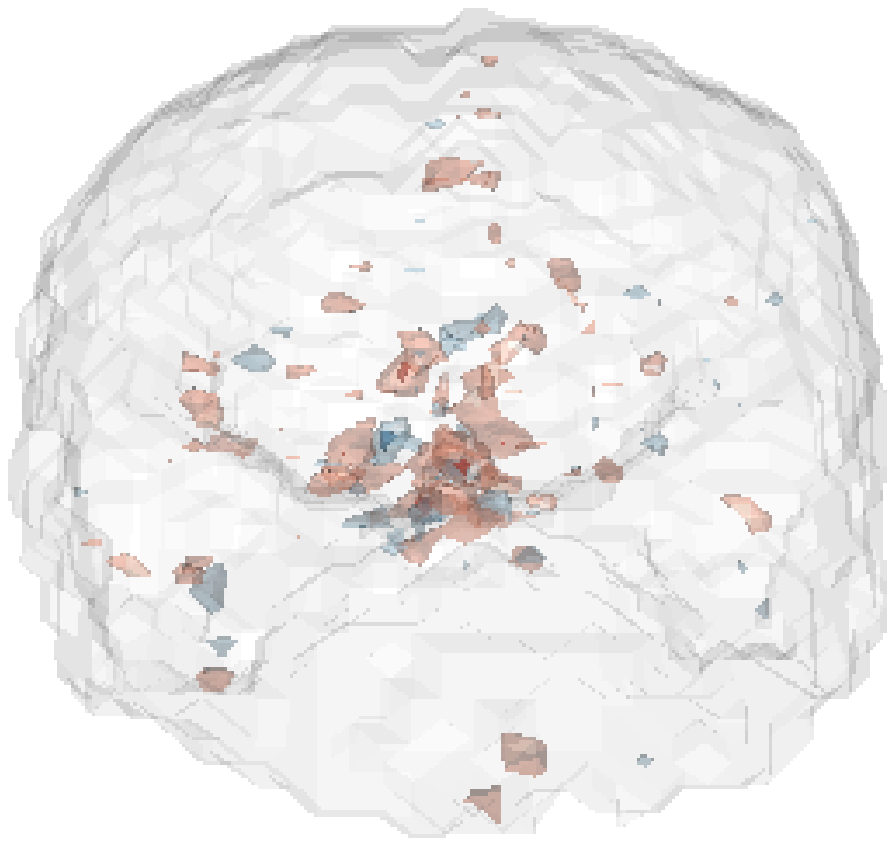} &
    \includegraphics[width=1.4in]{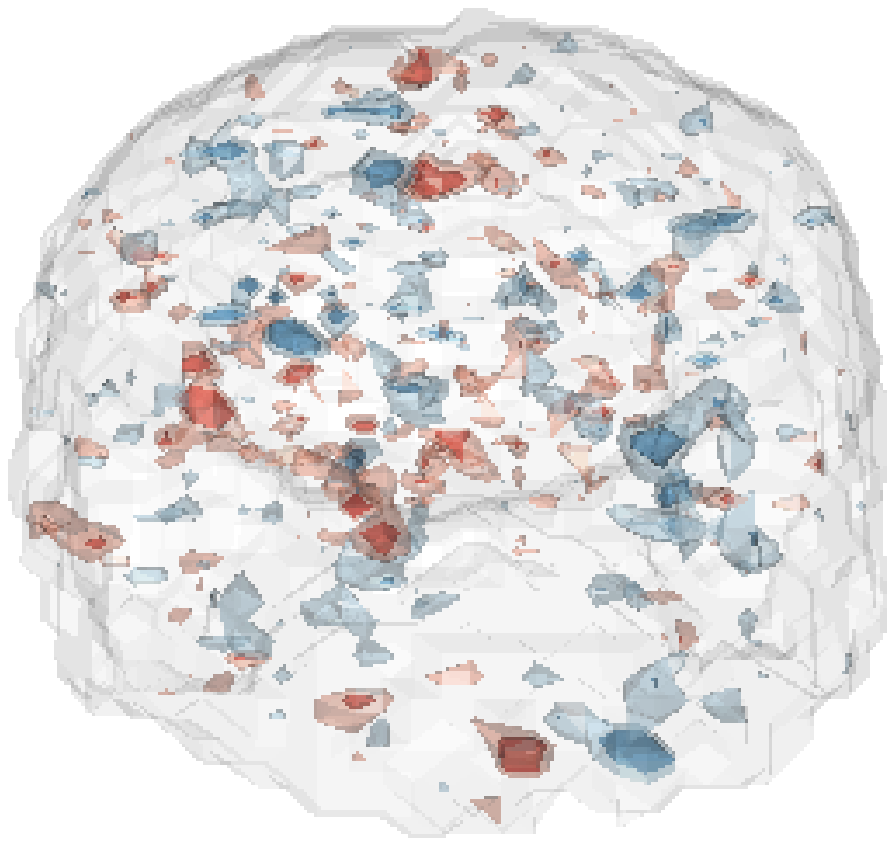} \\ 
    &  & 
  \end{tabular}
  \caption{Example of fMRI brain scans from the functional image analysis contest (FIAC).  Responses were measured for 16 subjects after hearing the same sentence spoken twice under two conditions: both spoken by the same speaker (top left) and each by a different speaker (lower left).  The gray shell gives the outline of the brain, while the light and dark red domains represent the contours at the hemodynamic activity levels (12,000, 15,000).  The difference brain (upper middle) shows the pointwise difference of these two scans.  In red regions, activity under the different-speaker condition exceeds that under the same-speaker condition.  In blue regions, the converse is true.  The pointwise average of these differences across the 16 experimental subjects is less noisy (lower middle).  The residual brain (upper right) gives the difference between the first subject and the experimental average.  Finally, the $t$ brain (lower right) shows paired $t$ statistics, where the light and dark colors represent the nominal single $t$-test 95\% and 99\% thresholds of 2.13 and 2.95.}
  \label{FigureBrainExample}
\end{figure}

A concrete example that we will study in this paper is an experiment involving functional magnetic resonance imaging (fMRI), which was a language priming experiment carried out by \citet{Dehaene2006} that appeared in the functional image analysis contest (FIAC).  fMRI responses were measured twice for 16 subjects after each heard a sentence spoken under two different conditions: once with the same speaker both times, and once with different speakers.  After each repetition, hemodynamic activity was measured at every point (voxel) in a 64 $\times$ 64 $\times$ 30 grid that encompasses the brain for each subject.  The fMRI scan for the first subject is shown in Figure \ref{FigureBrainExample}: the light and dark red domains represent mid- and high-activation regions of the brain, respectively, under the same-speaker (upper left) and different-speaker conditions (lower left).

This study inspires two important research questions: Firstly, are there any significant differences at all?  Secondly, if there are, where do these significant differences lie?  A natural way to look for differences is to compare the two conditions at every voxel in the brain.  Figure \ref{FigureBrainExample} (lower right) shows the result of voxelwise paired $t$ tests: in red regions, different-speaker activity exceeds same-speaker activity; in blue regions, the opposite is true.  The light and dark colored regions show areas where the $t$ statistics exceed the nominal single-test 95\% and 99\% levels of 2.13 and 2.95, respectively.  The multiplicity of comparisons of 24,759 voxels in the brain, each with its own $t$ test, poses a challenge in establishing significance.  Using the 95\% threshold, 1,273 voxels indicate significant differences, which, as a proportion of 24,759, amounts to approximately 5\%, and therefore mainly false positives.  In contrast to these excessive false positives using the $t$ test 95\% threshold, using the 95\% Bonferroni bound of 7.23 indicates there are no significant differences at all.  This bound is known to be very conservative, and the resulting conclusion of no significant results whatsoever is unconvincing, and moreover, unlikely to be inconsistent with neurological hypotheses.  Furthermore, in considering the nature of fMRI studies, the strength of spatial dependence depends on the resolution, among other factors, and should thus also be taken into consideration.  These preliminary exploratory results underline the importance of our motivating research questions and provide the impetus for a sound method to determine and characterize significance.

\subsection{Tail probabilities and Euler characteristics}
As a starting point, let the null hypothesis $H_{0}$ assert that both brain scans, under both conditions, are, on average, equivalent.  For the moment, assume that rather than having voxel-based data, we observe the fMRI images over a continuum.  Issues of resolution and discretization will be addressed and formalized later.  Under such an $H_{0}$, the grid of $t$-statistics is a smooth random field $T$ over a region $S$ of the brain.  Since high values of $T$ usually indicate deviation from $H_{0}$, a natural test statistic to consider is then the maximum
\begin{equation*}
  M_S \overset{\Delta }{=}\sup_{s\in S} T(s).
\end{equation*}
The problem with using $M_S$ as a test statistic, however, is that its null distribution is required, but is virtually never known.  In the FIAC data, $M_S$ = 6.08; a test of $H_0$ requires the $p$-value, $P(M_S \geq 6.08 | H_0)$.  In addition, to identify the activated regions of $T$, we require a 5\% threshold $t$ such that $P(M_S > t | H_0)$ = 0.05.  Obtaining these unknown quantities, the null distribution and the threshold, is difficult because they both depend on the correlation structure of $T$, which itself is also unknown.


To bypass this problem, we make use of a powerful parametric method for determining the null tail probabilities of $M_S$ for a wide class of Gaussian and Gaussian-related random fields: the Euler characteristic heuristic (ECH) \citep{Adler:2000}, which provides an accurate approximation to the exceedence probability $P(M_S \geq u)$ for large $u$ ({\em i.e.}~high levels or thresholds $u$), where ``large'' refers to values of $P(M_S \geq u)$ of the order of 5\% or smaller, and ``accurate'' refers to an error in approximation of the order of 1\% to 2\% in the threshold level. (i.e. A true threshold of, say, 5, in standardised units, will be estimated by a value in the range (4.99,5.01).)  The ECH has been widely used in topological inference in fMRI studies \citep[{\em e.g.},][]{Cao:Worsley:1999:Correlation,FRISTON94,Friston-Kilner,Shafie:Sigal:Siegmund:Worsley:2001,Worsley:1994:Chi:t:F,Worsley:1995:Boundary,WOR95a}, cosmological data \citep[{\em e.g.},][]{Bardeen,Sloan2008,gott2007genus,TORRES,VPGHG}, as well as a number of other areas.  Moreover, for smooth Gaussian random fields, the ECH is not merely a heuristic: it is a rigorous theory \citep*{Taylor:Takemura:Adler:2003:Validity, RFG}.

\begin{figure}[!htbp]
  \begin{tabular}{cccc}
    $\boldsymbol{\varphi(A_{-4}) = 1}$ &
    $\boldsymbol{\varphi(A_{0}) = -493}$ &
    $\boldsymbol{\varphi(A_{2}) = 98}$ &
    $\boldsymbol{\varphi(A_{3}) = 3}$ \\
    \includegraphics[width=1.1in]{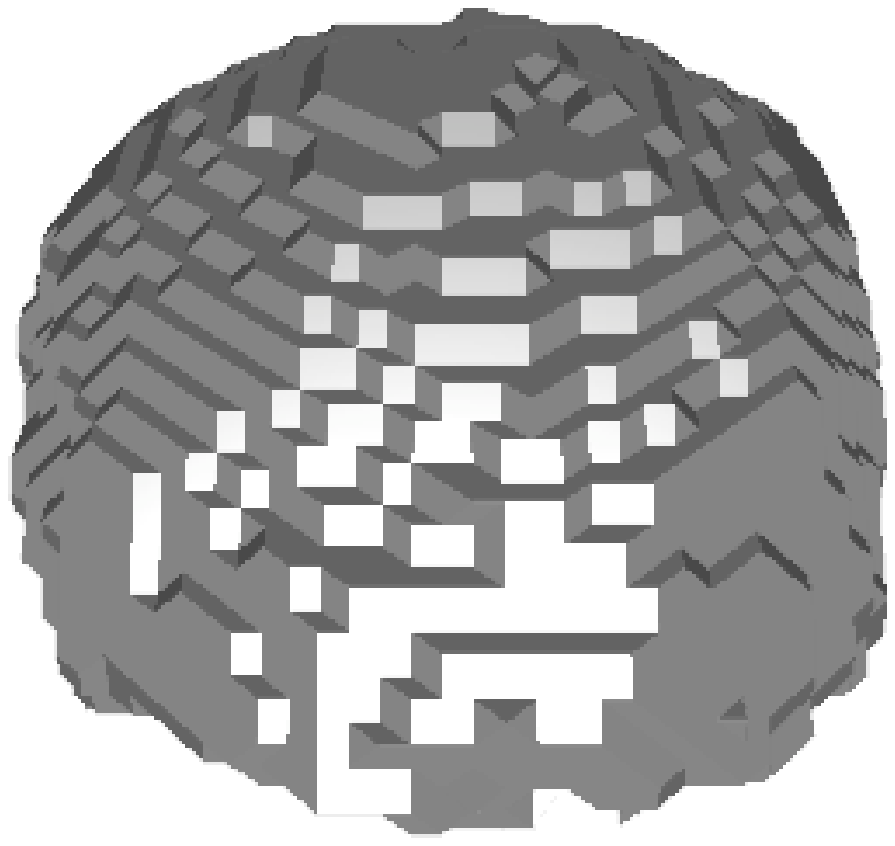} &
    \includegraphics[width=1.1in]{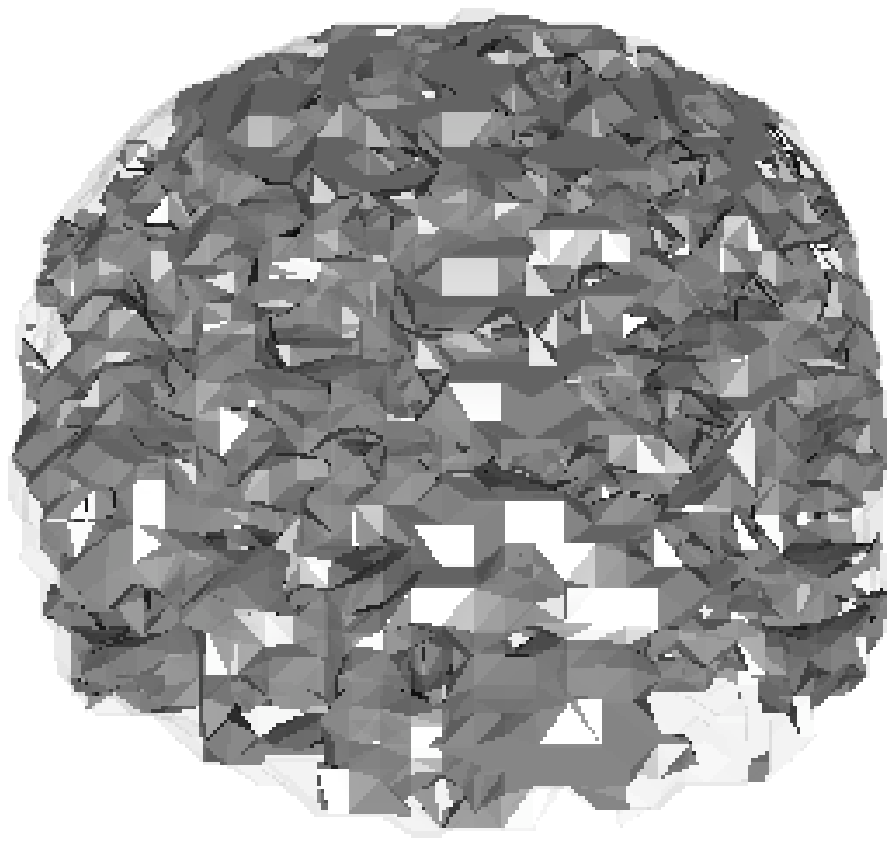} &
    \includegraphics[width=1.1in]{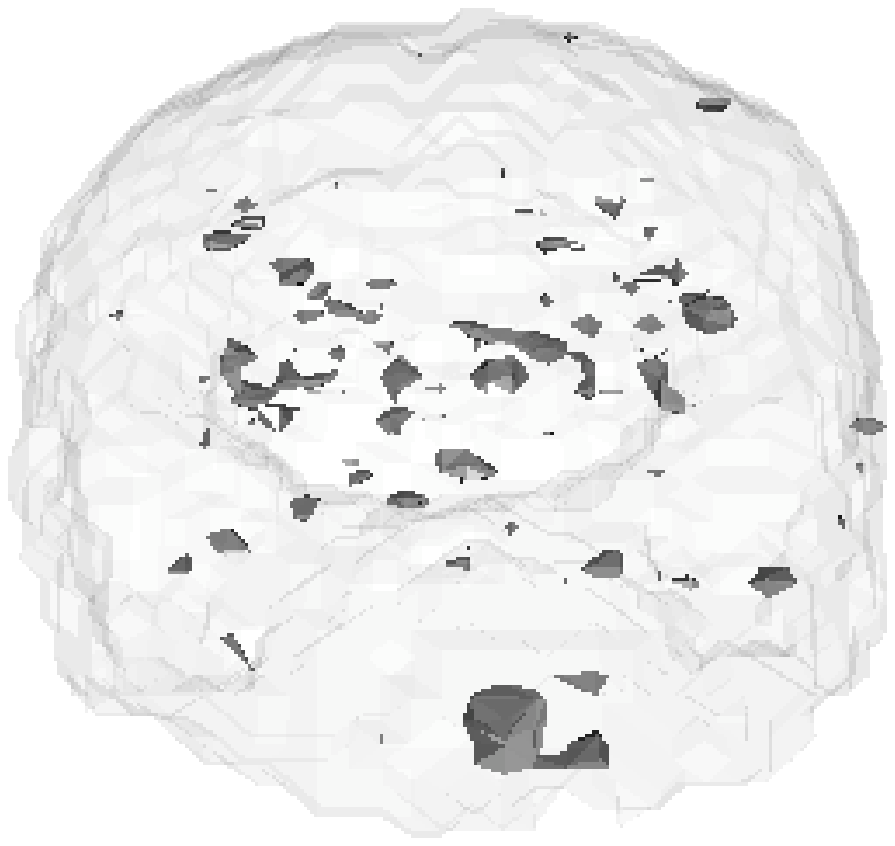} &
    \includegraphics[width=1.1in]{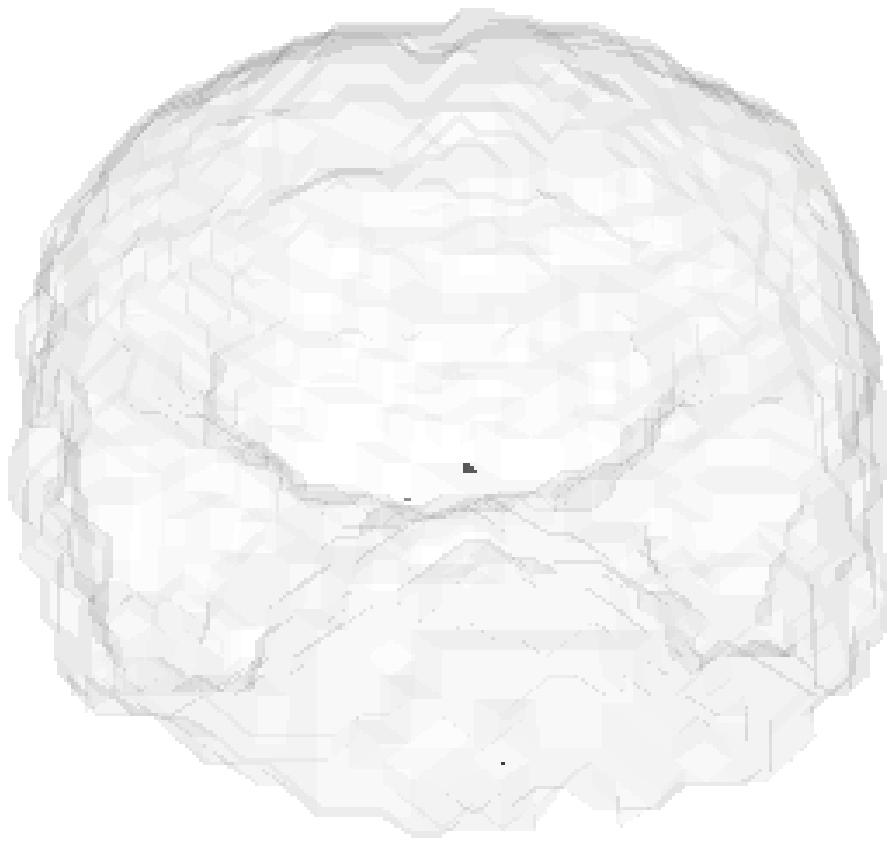}
  \end{tabular}
  \caption{Four excursion sets for a random field from the FIAC data.  Each excursion set is defined as $A_{u} \overset{\Delta }{=} \{s \in S : T(s) \geq u\}$.  The light gray shell shows the outline of the brain, while the dark gray regions show $A_u$ for $u$ = -4, 0, 2, and 3.}
  \label{FigureECExample}
\end{figure}

The ECH is based on the Euler characteristic (EC) $\varphi$, an important topological quantity for many general classes of well-behaved sets.  For a 3-dimensional Euclidean volume $V$, $\varphi$ counts the number of each of the three types of topological features of a manifold ---  (i) solid, simply connected regions of the manifold, or connected components, (ii) visible, open holes or handles,  and  (iii) invisible, closed holes or voids --- in an alternating sum.  It is a topological invariant, {\em i.e.}~a property of a topological space that is invariant under homeomorphisms, or continuous deformations, of the topological space.  It is given by 
\begin{eqnarray}
\label{EquationEC}
  \varphi(V) &=& \text{\# connected components in $V$} \\ &&\qquad - \text{(\# handles in $V$)}\  +\  \text{(\# voids in $V$)}. 
\notag
\end{eqnarray}
The ECH considers the ECs of a specific case, that of the excursion sets $A_u$,
\begin{equation}
  A_{u} \overset{\Delta }{=} \{s \in S : T(s) \geq u\};
  \label{EquationExcursionSet}
\end{equation}
Figure \ref{FigureECExample} shows examples of $A_u$ and $\varphi(A_u)$ for varying $u$; $\varphi(A_u)$ is an integer that can be negative (middle left), positive (middle right) or close to zero (rightmost).  Since topological features of the manifold lying below the level $u$ are effectively ignored, the manifold (and thus also the EC) becomes simpler as $u$ moves higher: At moderately high levels of $u$, $A_u$ typically takes on the simplified form of a union of simply connected components, while at the highest levels, it either disappears completely ({\em i.e.}~it is empty) or all that remains is a single simple component.  Thus, at such high levels, the EC is either zero or one.  Based on this phenomenon, the ECH claims that the expected EC for high $u$ approximates the tail probability,
\begin{equation}
  P(M_S \geq u) \approx E[\varphi(A_u)]. \label{EquationMax}
\end{equation}
A further result based on this phenomenon of topological simplification of excursion sets with increasing $u$ is that the expected EC is guaranteed to lie in $[0,1]$ for $u$ high enough.  From (\ref{EquationMax}), we obtain the required $p$-value for $H_0$, which can then be inverted to find the $u_\alpha$, or the threshold, corresponding to a $(1-\alpha)$ confidence level, via the approximation
\begin{equation}
  \widehat{u_\alpha}\overset{\Delta }{=}
  \max \{u : E[\varphi(A_u)] \geq \alpha\}.
  \label{EquationThreshold}
\end{equation}

When considering random fields \citep[see][]{RFG} with constant mean and variance (which for the moment we take to be zero and one), the ECH provides the additional practicality of a parametric closed form for the expected EC $E[\varphi (A_{u})]$, given by the Gaussian kinematic formula (GKF) \citep{Taylor2006},
\begin{equation}
  E[\varphi (A_{u})]=\sum_{i=0}^{\text{dim}(S)}\mathcal{L}_{i}(S)\rho _{i}(u).
  \label{EquationGKF}
\end{equation}
Here the $\rho_i(u)$ are functions that usually take on a simple and explicit form. For example, if the random field is Gaussian, 
then the $\rho_i(u)$ are explicitly expressible in terms of Hermite polynomials, as we shall see in 
\eqref{rho-rob} below.  Other examples can be found in \citet{GRF,ATSF,ARF}. The functions $\mathcal{L}_{i}(S)$, which depend on both the domain $S$ and the covariance of $T$, are the Lipschitz-Killing curvatures (LKCs), or intrinsic volumes, of $S$ --- complex topological quantities that are often extremely difficult to evaluate theoretically, even for dedicated topologists, as well as numerically.  When known, however, the LKCs $\mathcal{L}_{i}(S)$ provide the expected EC for all $u$ via the GKF (\ref{EquationGKF}), which can in turn be substituted into (\ref{EquationMax}) provided by the ECH to obtain an approximate $p$-value for $H_0$ or into (\ref{EquationThreshold}) to obtain a threshold level.  The estimation of the LKCs  is the bottleneck of a successful application of the ECH and is the central concern of this paper.

\subsection{Estimating \LK\ curvatures}
\begin{figure}[!ht]
\includegraphics[width=2.0in]{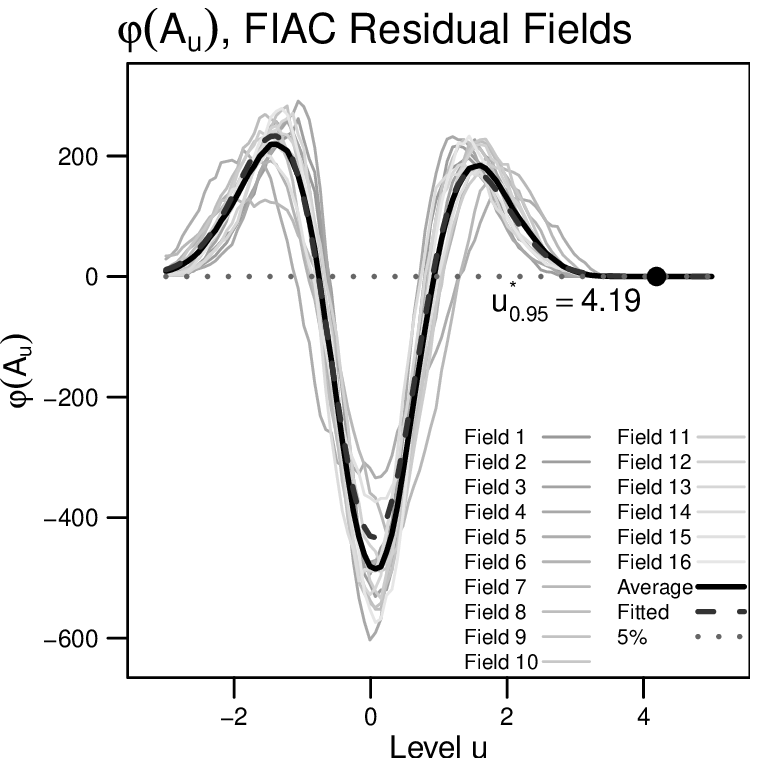} \hspace*{0.2in}
\includegraphics[width=2.0in]{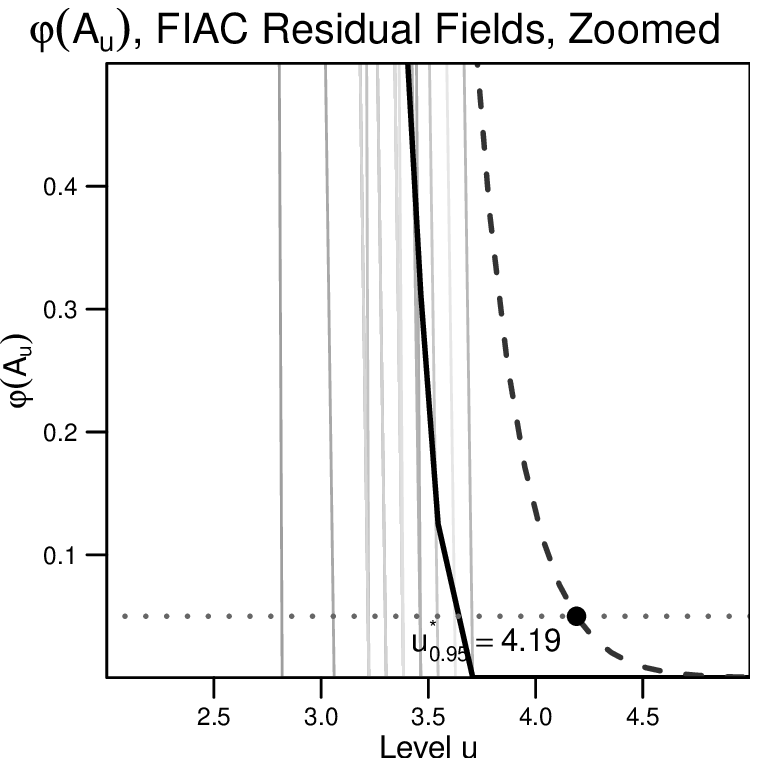}
\caption{Illustration of the proposed regression method to fitting Lipschitz-Killing curvatures (LKCs).  The light grey lines show the observed Euler characteristic (EC) profiles for each of the 16 observed fields from the FIAC data.  The solid black line is their average.  The dashed black line is the expected EC estimated by our regression method.  The black dot is the point where the expected EC intersects with dotted grey line at 0.05, which denotes the 5\% threshold $u_{95\%}^{\ast }$.  The plots are displayed both over a broad range of levels (left) and zoomed in near the threshold (right).}
\label{FigureBrainLKC}
\end{figure}

We introduce a new approach that estimates the LKCs by simply matching the expected and empirical ECs, which then directly yields the $p$-value and threshold level in a straightforward manner via (\ref{EquationMax}), (\ref{EquationThreshold}), and (\ref{EquationGKF}).  Figure \ref{FigureBrainLKC} provides an illustration of our method: the left panel shows the empirical $\varphi(A_{u}^{(i)})$ profiles for the 16 observed fields from the FIAC data (thin gray lines) and their sample average (solid black line) for different values of $u$.  We find the best-fitting LKCs through a generalized least squares regression of (\ref{EquationGKF}).  The fitted LKCs then produce the best-fitting profile: the dashed black line.  The intersection of this fitted profile and the 0.05 line yields the 95\% confidence threshold as $u_{95\%}^{\ast}$ = 4.19 (black dot).  This threshold can be applied to the brain of $t$ statistics to identify significantly activated regions.  Doing this indicates only 7 voxels activated beyond 4.19 --- a striking contrast to the 1,273 found using the naive  $t$-test 95\% threshold of 2.13, and also more credible than the 0 found using the 95\% Bonferroni bound of 7.23.

One caveat of which to take note in matching expected and empirical ECs is that the observed (empirical) ECs $\varphi (A_{u}^{(i)})$ at very high levels of $u$ are generally too noisy to directly estimate the threshold $u_{95\%}^{\ast}$; calculating a threshold from extreme empirical ECs is akin to estimating a tail quantile with a small number of observations (16, for the FIAC data).  The right panel of Figure \ref{FigureBrainLKC} illustrates this phenomenon by zooming in and showing that the empirical and expected ECs cross the 0.05 line at very different values of $u$.  In cases of very high $u$, for general $p$-value approximation and threshold level determination, use of the expected EC yields significantly more reliable and accurate results that do not depend on small sample sizes, but instead derive from rigorous probability theory and topology. 

Inference via LKC regression is parametric because it relies on the ECH, which leverages on an assumed underlying Gaussianity.  A commonly suggested nonparametric alternative is a permutation test on the group assignments (same-speaker, different-speaker) to approximate the null distribution of $M_S$.  The main drawback of this nonparametric option is the tremendous expense of computing the maximum $t$-statistic given a set of permuted labels, which is severely impeded by the computation of many hundred thousand voxelwise averages and standard deviations.  For the FIAC data, a permutation test with just 100 permutation samples takes several hundred times longer than LKC regression.  Moreover, observed random fields are typically Gaussian-related by construction, since they arise as the residual of some statistical procedure.  Nonparametric methods generally underperform in cases where such distributional properties are known.

To date, the only other existing estimation method for LKCs is the warping method of \citet{Taylor2007}.  Warping transforms the realizations $T^{(1)}, T^{(2)}, \ldots, T^{(F)}$ of the random field to a volume in a $F$-dimensional space, where $F$ is the number of random fields in the data; for the FIAC data $F$=16.  The intrinsic volumes of the transformed volume are shown to be good estimates of the LKCs.  Though accurate, understanding the method requires a knowledge of Riemannian geometry, and can be computationally slow, particularly for high dimensions.  Compared to \cite{Taylor2007}, our method has similar accuracy but is computationally faster.  Furthermore, it requires only the ECs of the observed random fields.  Since observations are typically taken over a regular grid, efficient numerical optimization is available.  The warping method, on the other hand, must perform the slower computation of intrinsic volumes on a high-dimensional transformed space with irregularly spaced points.

\subsection{Paper structure and acknowledgements}
The remainder of this paper is organized as follows: In Section \ref{SectionLKCRegression}, we introduce our regression-based estimation method, detailing on the various implementation factors and combinations available, and their associated issues of optimality.  Section \ref{SectionSimulationExperiment} presents the results of a large scale simulation study, which compares these various factors and their combinations, as well as our approach to the warping method of \citet{Taylor2007}.  In Section \ref{SectionApplication} we extend this comparison from simulation studies to the concrete example of the FIAC data.  Section \ref{SectionConclusion} concludes the paper with a summary of its main results and proposals for other areas where our method would be applicable.  The Appendix provides a brief treatment of the Gaussian kinematic formula (\ref{EquationGKF}).

Finally, we would like to thank Jonathan Taylor for helpful discussions at various stages of this work, and Naor Alush for help with the programming.


\section{Lipschitz-Killing Curvature Regression}

\label{SectionLKCRegression}

The goal of our regression method is to obtain accurate estimates of the Lipschitz-Killing curvatures (LKCs) $\mathcal{L}_{i}$ for use in $p$-value calculation and threshold level determination (see (\ref{EquationExcursionSet}) and (\ref{EquationThreshold})).  (While the LKCs are of interest in and of themselves, in this paper we shall concentrate almost exclusively  on their  importance for threshold determination.) 
The inputs are the observed realizations $T^{(i)}$, $i=1,\ldots ,F$, of a random field, assumed to be Gaussian or Gaussian-related and normalized to mean zero and unit variance, over a region $S$.

The Gaussian kinematic formula (GKF) (\ref{EquationGKF}), combined with the central limit theorem, suggests a linear model with heteroscedastic, possibly correlated errors: 
\begin{eqnarray}
\frac{1}{F}\sum_{i=1}^{F}\varphi (A_{u}^{(i)}) &=&\sum_{i=0}^{\dim (S)}%
\mathcal{L}_{i}\rho _{i}(u)+\varepsilon (u),  \label{EquationModel} \\
\varepsilon (u) &\sim &N(0, \sigma_u^2),  \notag
\end{eqnarray}
for some unknown $\sigma_u^2 = \text{Var}(\varphi\big(A_u)\big)/F$.  The response variable is the average of the empirical Euler characteristics (ECs).  The regressors $\rho_i(u)$ take on a simple and explicit form; when the random field is Gaussian,
\begin{equation}
\label{rho-rob}
\rho_i(u)=(2\pi)^{-(i+1)/2}H_{i-1}(u)e^{-u^{2}/2},
\end{equation} 
where $H_{j}$ is the $j$-th Hermite polynomial.  We estimate the unknown LKCs $\mathcal{L}_{i}$ through generalized least squares (GLS).

Once the LKCs are estimated, we easily obtain an estimate of the expected EC via the GKF (\ref{EquationGKF}) for any level $u$,
\begin{equation*}
  \hat{E}[\varphi (A_{u})]=\sum_{i=0}^{\text{dim}(S)}\mathcal{\hat{L}}%
  _{i}(S)\rho _{i}(u),
\end{equation*}
which then substitutes into (\ref{EquationMax}) to give an approximate $p$-value for hypothesis testing for large $u$, and into (\ref{EquationThreshold}) to give a
threshold level.

The validity of the GKF (\ref{EquationGKF}) for all $u$ forms the foundation of the effectiveness of our LKC regression method.  For $p$-value and threshold calculation, one typically encounters large values of $u$, which correspond to small-probability events, the direct estimation of which is unreliable.  Our regression approach, instead, is grounded on the observation that when $u$ is small or moderate, the expected EC $E[\varphi (A_{u})]$ can be well estimated from the data, since these cases do not correspond to small-probability events.  These reliable estimates then translate, in turn, through our regression, into reliable estimates of the LKCs $\mathcal{L}_{i}$, which do not depend on $u$; the $\mathcal{\hat{L}}_{i}$ then yield good approximations for $p$-values and threshold levels for large $u$.  In summary, LKC regression leverages the precision of estimation at low levels of $u$ to obtain accurate approximation at high levels of $u$.

Several issues need to be addressed in the construction of our LKC regression method: First, to evaluate the empirical EC, the region $A_{u}^{(i)}$ for fixed $u$ must be determined, which poses a potential difficulty because each field $T^{(i)}$ is typically observed at a discrete grid of sample sites, yet, to apply the GKF, the excursion sets $A_{u}^{(i)}$ must be described as smooth regions over $S$.  We detail the construction of $A_{u}^{(i)}$ in Section \ref{SectionObservedExcursionSets}.  Second, compared to conventional linear regression, an interesting feature of our setting is that arbitrary amounts of ``data'' can be generated at very little cost, simply by exploiting (\ref{EquationModel}) at a collection of specified levels $u$.  Each $u$ gives distinct excursion sets $A_{u}^{(i)}$ and their corresponding average EC, 
\[ \frac{1}{F}\sum_{i=1}^{F}\varphi(A_{u}^{(i)}).
\]
However, this supposed advantage in fact raises the problem of a tradeoff between residual error in the regression and estimation error in the covariance matrix: additional data might help the regression, but the covariance matrix becomes increasingly difficult to estimate with more data.  We address this tradeoff in detail in Section \ref{SectionDesignSelection}.  Third, after specifying the levels $u$, the unknown error covariance matrix must be estimated before running GLS.  We present several estimation options in Section \ref{SectionErrorCovariance}.  In Section \ref{SectionSimulationExperiment} we find optimal settings for level selection of $u$ and error covariance estimation.

As a tool for hypothesis testing, the LKC regression approach is powerful because it overcomes the need for the complete specification of the covariance structure of the random field.  Very few assumptions about covariance are made; even isotropy is not required.  The random fields only need be Gaussian or Gaussian-related, which is often a natural consequence of the data generating process.

\subsection{Observed Excursion Sets}

\label{SectionObservedExcursionSets}

\begin{figure}[!htbp]
\includegraphics[width=1.2in]{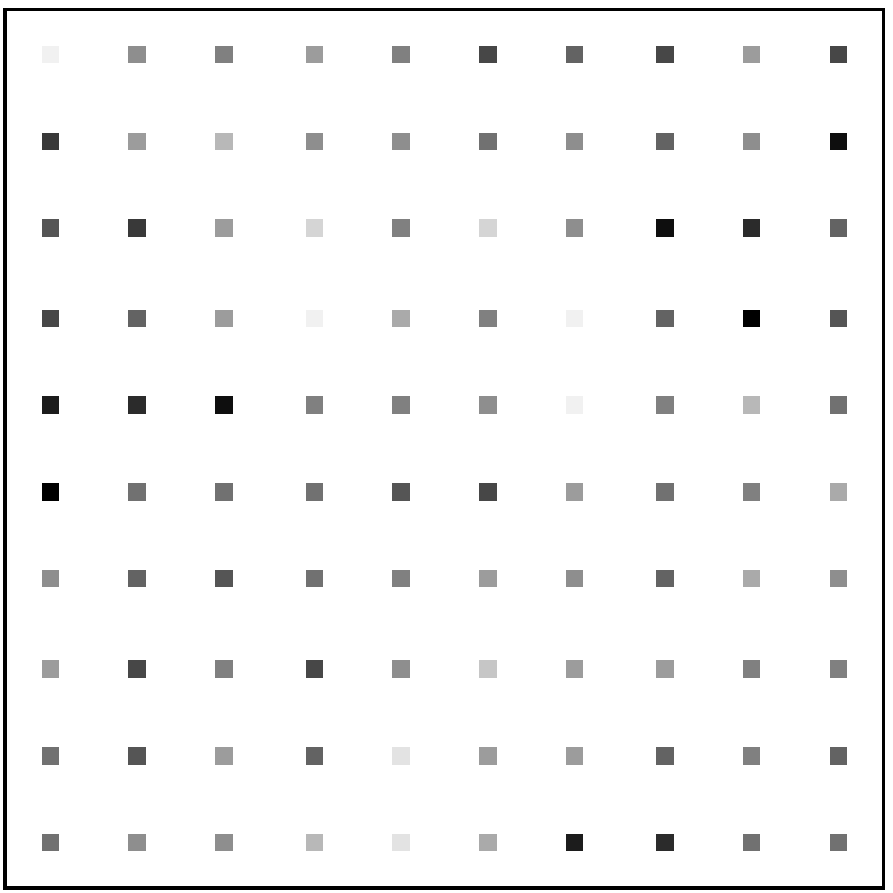}  \hspace*{0.2in}
\includegraphics[width=1.2in]{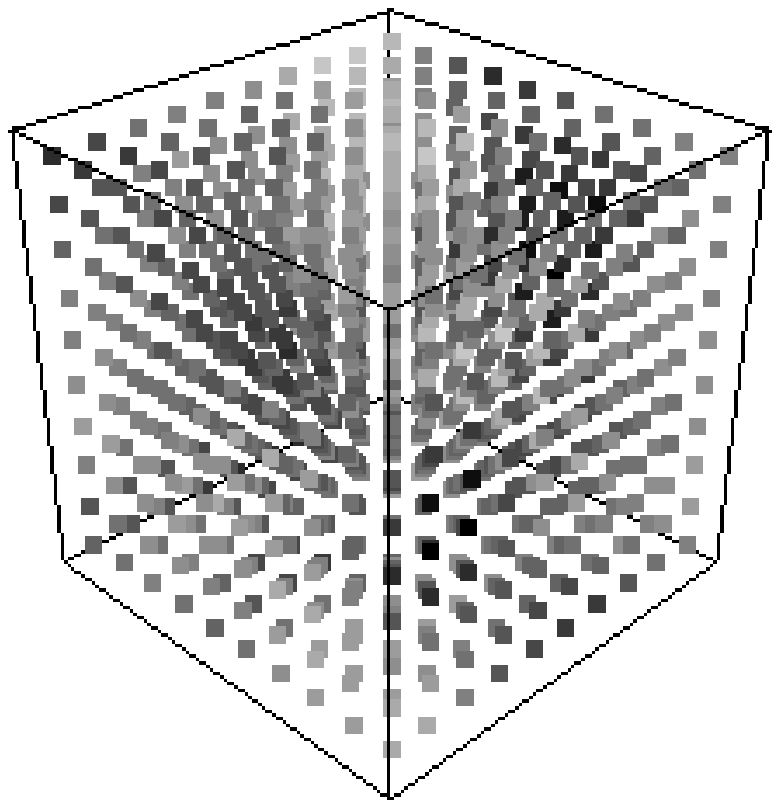}  \hspace*{0.2in}
\includegraphics[width=1.2in]{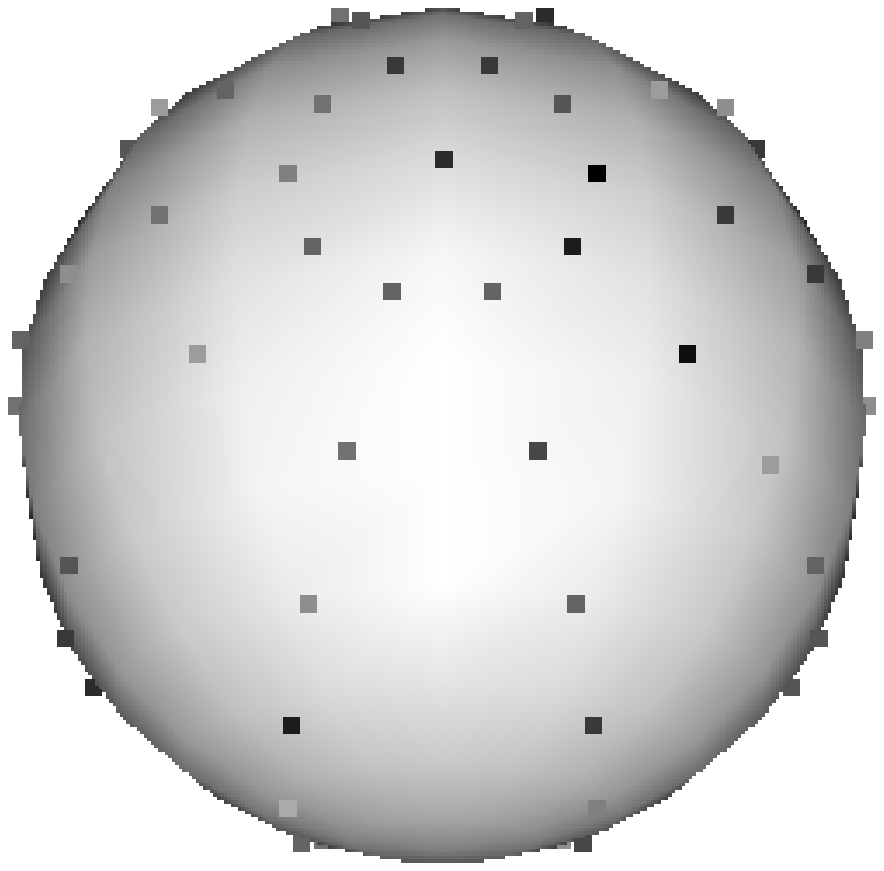} \\ \vspace*{0.2in}
\includegraphics[width=1.2in]{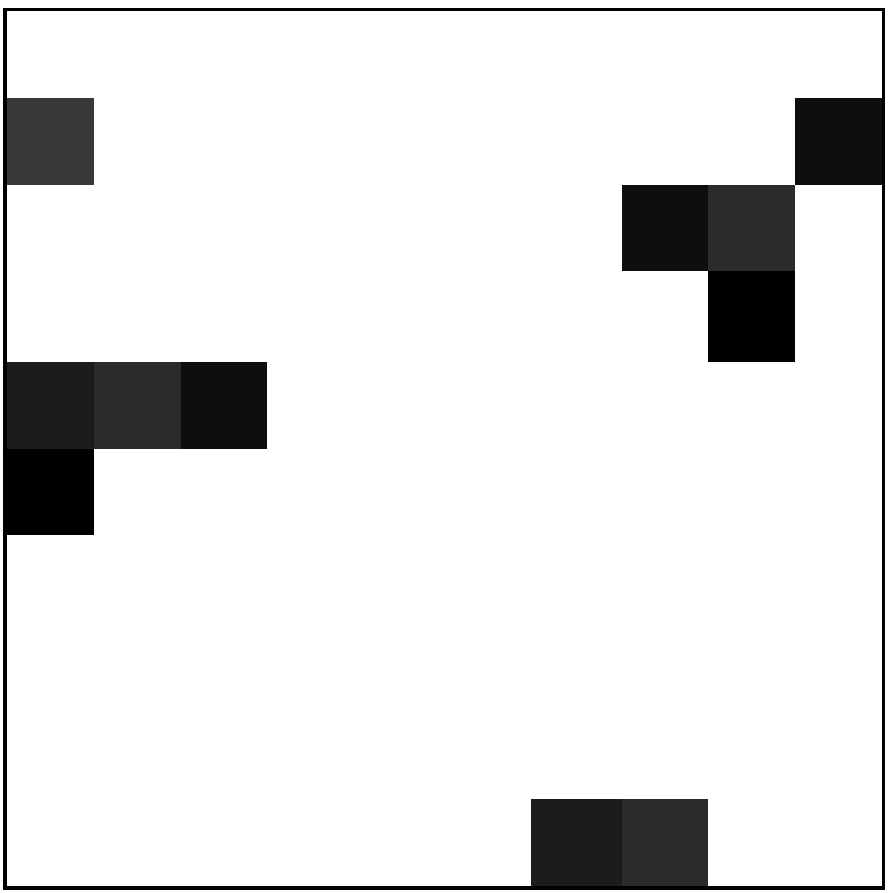}  \hspace*{0.2in}
\includegraphics[width=1.2in]{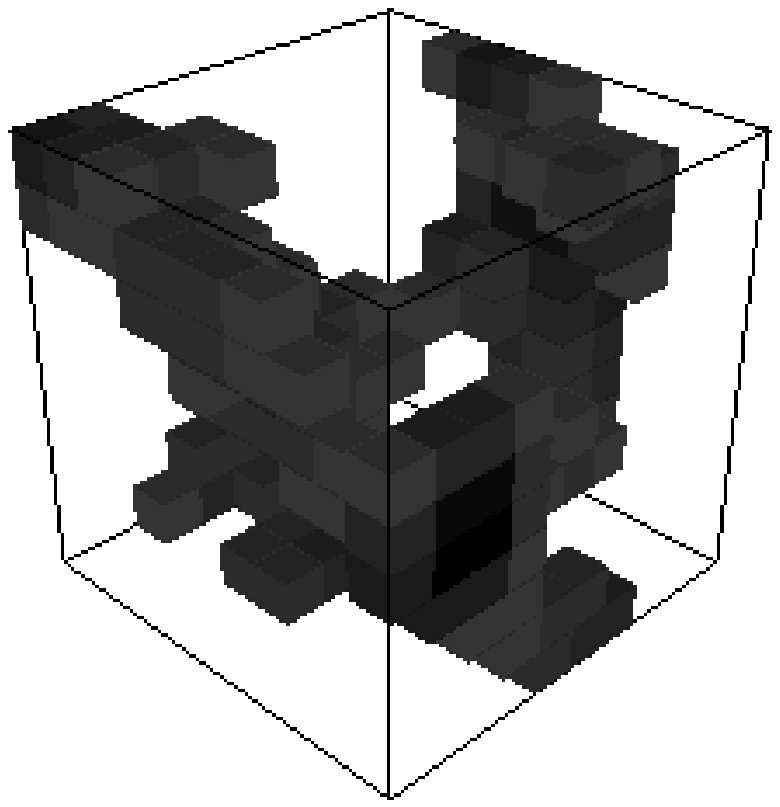}   \hspace*{0.2in}
\includegraphics[width=1.2in]{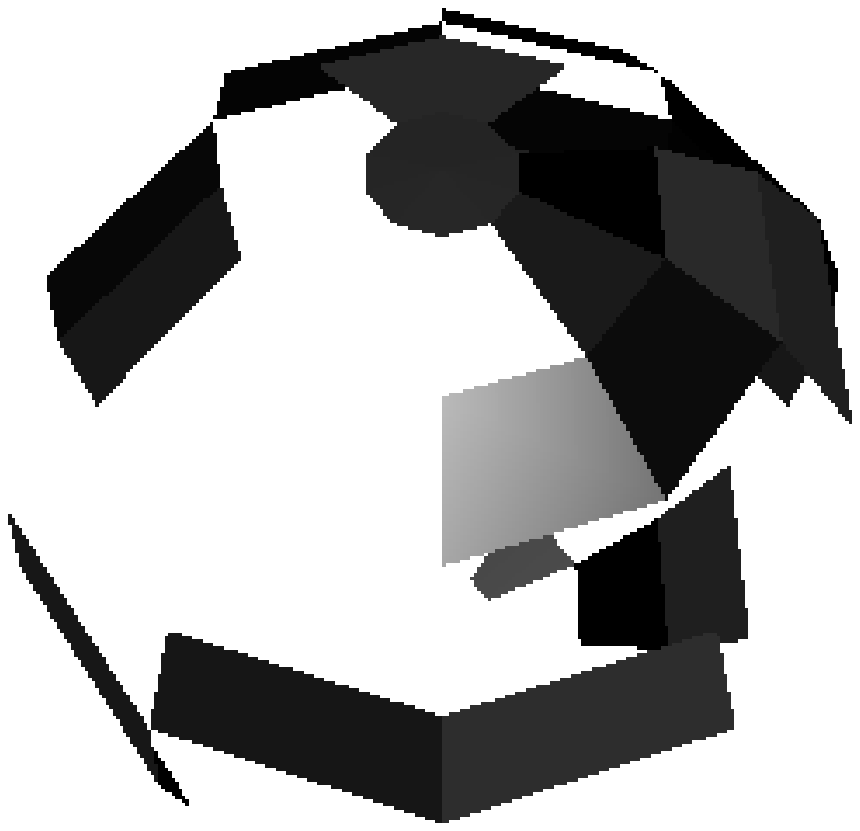}
\caption{Simulated examples of observed fields on the unit square (left), on the unit cube (middle) and on the unit sphere (right).  The two rows show the sample sites (upper row) and approximated excursion sets (lower row).  In each case darker points denote greater values.  The simulated random fields have Gaussian covariance and are simulated over a discrete grid with $G =$ 10 sample sites in each dimension for the square and cube and $G =$ 50 for the sphere.}
\label{FigureObservedExample}
\end{figure}

Though assumed to be continuous, in practice random fields are typically observed at a discrete set of sample sites $s_{k}\in S,\ k=1,\ldots ,K$.  An empirical realization of a random field $T$ is thus defined by the observations $\{T(s_{k}),k=1,\ldots ,K\}$; Figure \ref{FigureObservedExample} provides an illustration.  Sample sites commonly lie on a rectangular lattice for fields on the square and cube (upper left and middle) or on a latitudinal-longitudinal lattice for fields on the sphere (upper right).  We measure the resolution of the field by the grid size $G$, the number of sample sites in each direction/dimension; for example, for a 3D cube $G=$20 corresponds to $K=G^{3}=$8000.

In theory, to compute the EC $\varphi$ of an excursion set $A_{u}^{(i)}$, the excursion set needs to be represented as a smooth domain.  However, as the fields are observed only at discrete points $s_{k}$, an approximation of the exact domain is required.  In the digital topology literature, there exist many constructions \citep[e.g.,][]{Imiya1999} of domains created from a set of points, known as quasi-objects.  The methods are distinguished by their connectivity, which specifies how close two sample points must be for them to be considered part of the same domain.  For example, in considering points on a planar lattice, 4-connectivity places vertically and horizontally adjacent points into the same domain, but diagonally adjacent points into separate domains (unless they also share an immediate neighbor); 8-connectivity includes all diagonal adjacencies in the same domain.

To create excursion sets from observed random fields, we use planar (2D) 8-connectivity and spatial (3D) 26-connectivity, where each point touches all its immediate and diagonal neighbors.  This is equivalent to building $A_{u}^{(i)}$ by introducing a small square or cube centered at $s_{k}$ whenever $T^{(i)}(s_{k})>u$, and then joining them to form the excursion set.  For spherical fields, we use 8-connectivity over the latitudinal/longitudinal grid of sample sites, with the excursion sets made up of small plates.  Figure
\ref{FigureObservedExample} shows examples of the resulting excursion sets for fields on the plane (lower left), cube (lower middle) and sphere (lower right).  We experimented with different connectivity structures, and concluded that they had little impact on the calculated ECs.  In the brain image analysis example (FIAC, 64 $\times$ 64 $\times$ 30), we observed nearly identical EC profiles for spatial 6-connectivity, 18-connectivity, and 26-connectivity.

The final step in computing the response variable in the regression of (\ref{EquationModel}) is to evaluate the EC for excursion sets.  (\ref{EquationEC}) provides a simple formula, provided we can count the number of connected components, handles, and voids.  This calculation is easy in the discrete setting, in which the excursion set is made up of a large number of simple cells, such as squares, cubes or plates; the formula thus reduces to an alternating sum over distinct vertices, edges, faces and cubes of these cells:
\begin{equation}
  \varphi = \text{\# vertices} - \text{\# edges} + \text{\# faces} - \text{\# cubes} \label{EquationECCube}
\end{equation}
The same formula also holds in two dimensions without the final term.  A key advantage of our LKC regression method is the rapidity of this computation for cells on a grid; the warping method of \citet{Taylor2007}, in contrast, must perform the slower computation of all intrinsic volumes (not just the EC) over an irregular (nonrectangular), warped grid.

For an application with $F$ realizations of the random field, LKC regression calls for the computation of ECs at $U$ different excursion levels.  Taking advantage of binary search and careful indexing, the computation cost of our method is only $O(FK\log U)$.  This makes it extremely fast to evaluate the EC for many levels $U$.  

\subsection{Covariance of the Error Terms}

\label{SectionErrorCovariance}

The covariance structure of the error terms in LKC regression is defined by 
\begin{equation*}
C_\varphi(u,v)\overset{\Delta }{=}\text{Cov}(\varphi (A_{u}),\varphi (A_{v})).
\end{equation*}%
We emphasize that $C_\varphi$ should not be mistaken for the covariance of the random fields $T^{(i)}$; $C_\varphi$ describes the correlation of the Euler characteristics across differing levels, and not of field values themselves.  $C_\varphi$ can be viewed as the covariance function of the stochastic process $\varphi (A_{u}^{(i)})$ in $u$.

\begin{figure}[!ht]
     \includegraphics[width=2.2in]{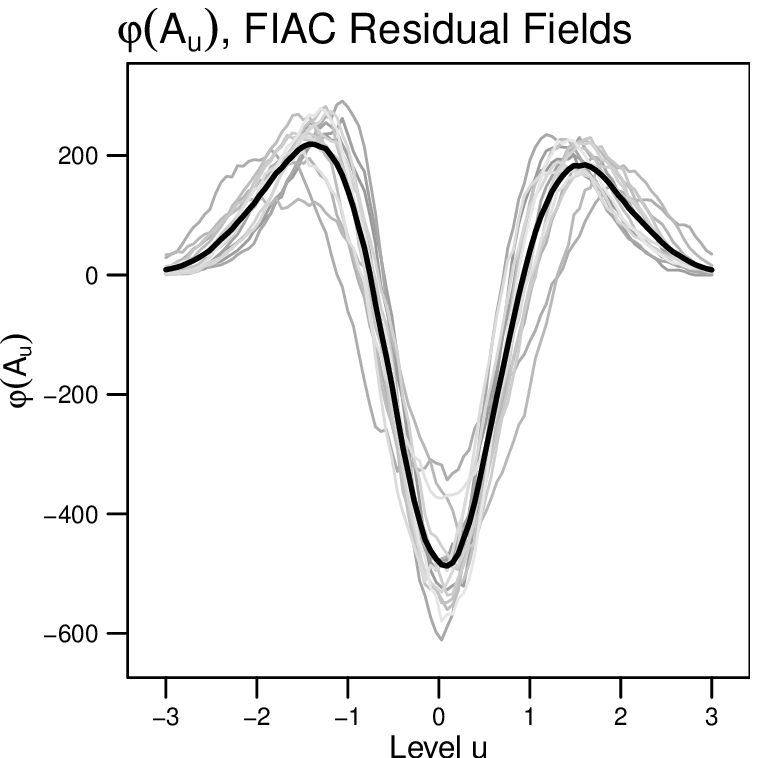}  \hspace*{0.2in}
         \includegraphics[width=2.2in]{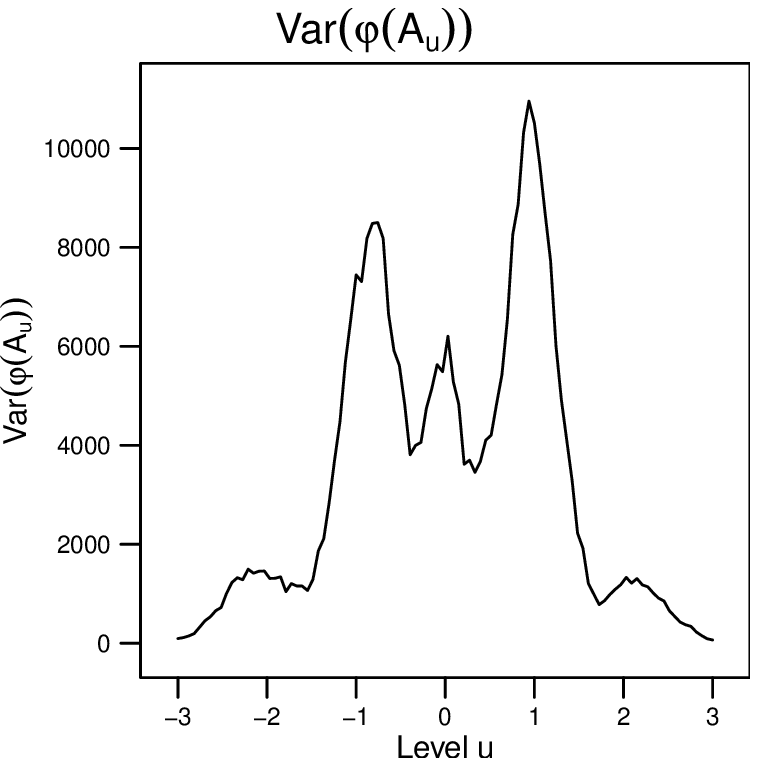} \\  \vspace*{0.2in}
    \includegraphics[width=2.2in]{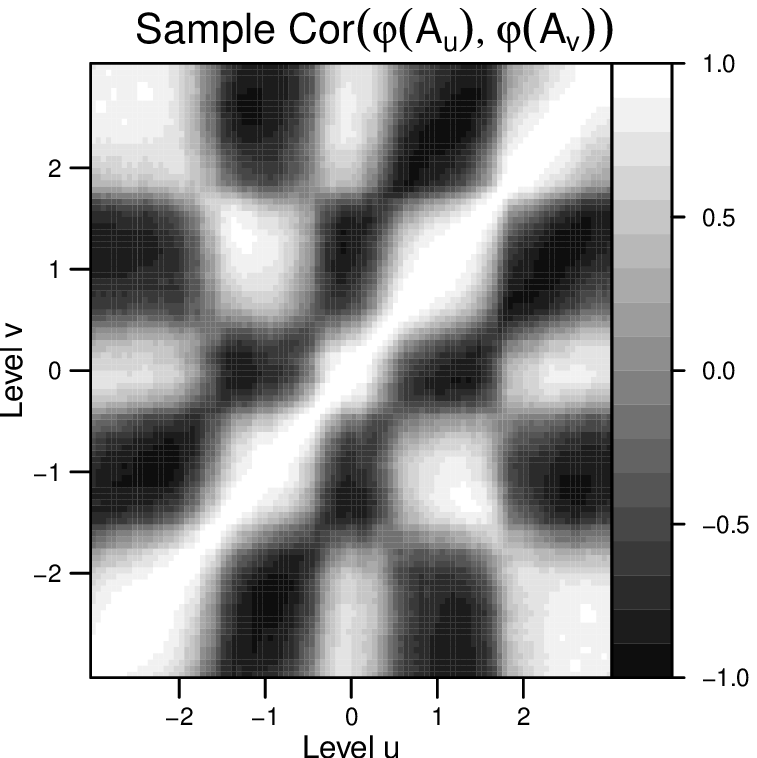}  \hspace*{0.2in}
        \includegraphics[width=2.2in]{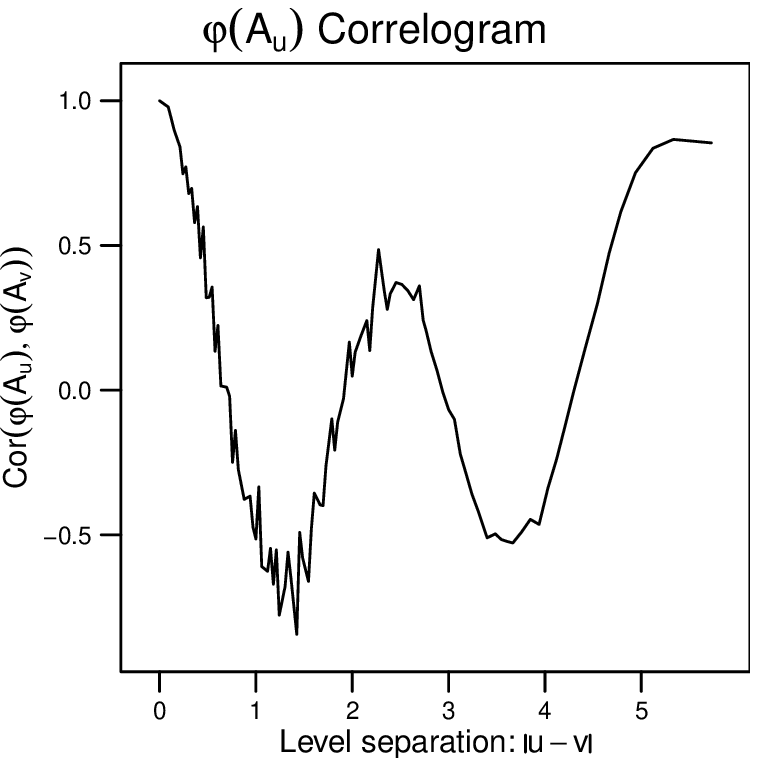}
  \caption{Illustrations of the profiles (upper left), variance (upper right), correlation (lower left) and correlogram (lower right) of the observed EC $\varphi(A_u)$ for the 16 FIAC fields.  In the correlation plots, the grayness reflects the strength of the correlation between two levels $u$ and $v$.  The correlogram reports average correlation as a function of $|u - v|$.}
  \label{FigureErrorCovariance}
\end{figure}

The error covariance structure of LKC regression is nontrivial: First, they are heteroskedastic; the observed EC tends to be less variable at higher levels of $u$, as illustrated in Figure \ref{FigureErrorCovariance} (upper right) for the FIAC data.  Second, they are correlated; the sample correlation plot (lower left) and correlogram (lower right) show pronounced negative correlation for levels separated by about 1.  This is due to the periodicity of the three swings in the EC profiles, which are about one (standardised)  unit apart.

\begin{figure}[!htbp]
  \includegraphics[width=1.5in]{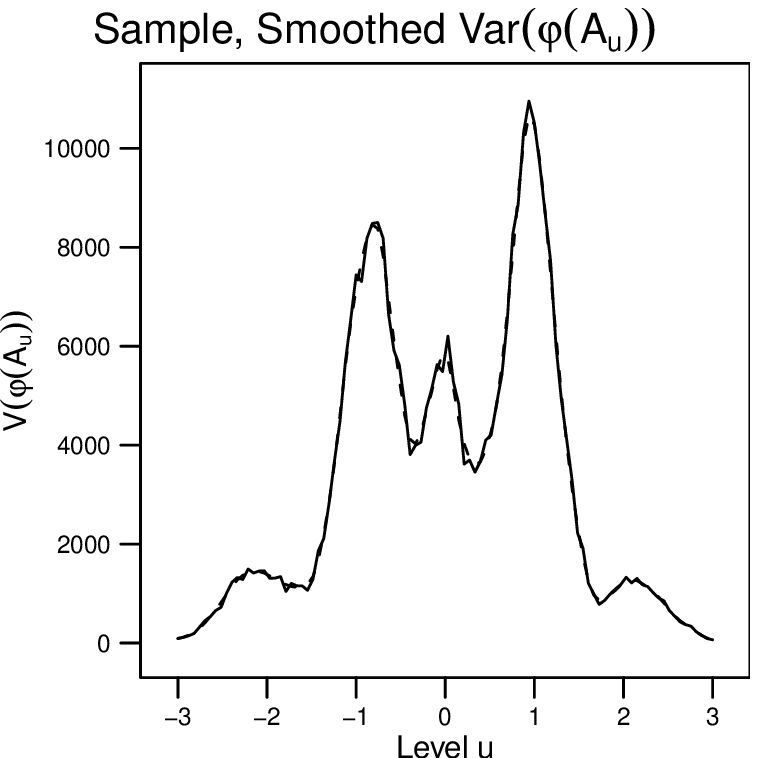}  \hspace*{0.2in}
    \includegraphics[width=1.5in]{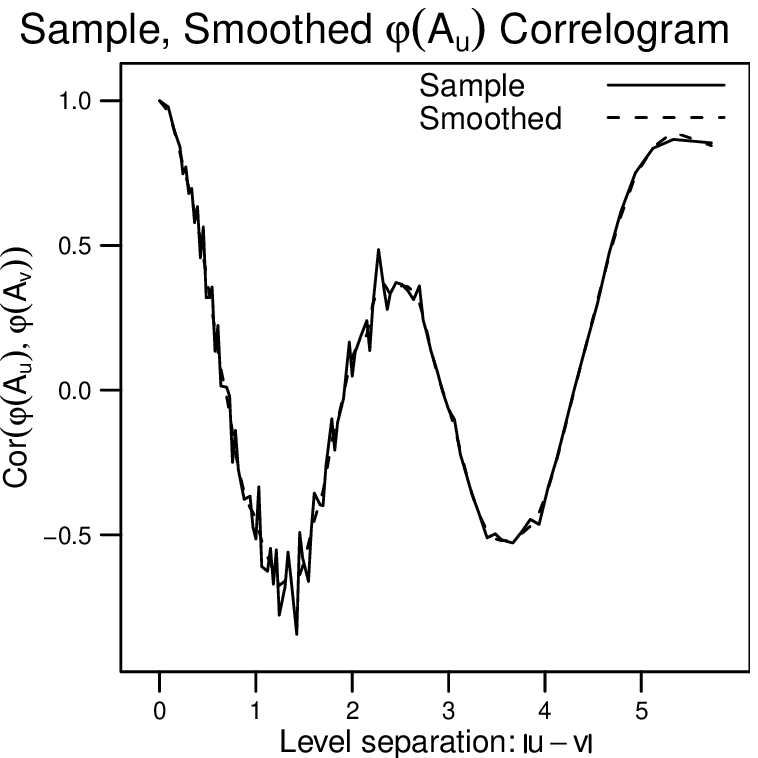}
  \caption{Smoothing applied to the variogram (left) and correlogram (right) by the smoothed diagonal (SD) and smoothed covariance (SC) methods.  Each plot shows both sample (dashed) and smoothed (solid) estimates of the variance and correlation of $\varphi(A_{u})$.}
  \label{FigureInternals}
\end{figure}

Given $U$ levels $u_{j}$, LKC regression requires the covariance matrix $\Sigma =[C_\varphi(u_{j_{1}},u_{j_{2}})]_{U\times U}$.  Although a misspecified covariance matrix still gives consistent estimates under GLS, the efficiency of the estimate depends on the estimation of $\Sigma$.  There are five ways of producing a positive-definite estimates of the covariance matrix $\hat{\Sigma}$:
\begin{enumerate}
\item {\it Identity (I).}  As a baseline, heteroskedasticity and correlation are ignored under this option so that $\hat{\Sigma}\propto I_{U}$, yielding an ordinary least squares (OLS) model.
\item {\it Smoothed diagonal (SD).}  Heteroskedasticity is incorporated while correlation is ignored.  The sample variances $\hat{V}_{u}=\widehat{\text{Var}}[\{\varphi (A_{u}^{(i)})\}_{i}]$ are computed and then smoothed as a function of $u$; we employ local quadratic smoothing with 10\% nearest neighbors \citep{Loader1999,Loader2010}.  To prevent the smoothing procedure from generating negative variances, the smoothing is applied to the logarithm of the sample variances; values of $u$ with a sample variance of zero are dropped from the smoothing procedure.  Examples of original and smoothed variances are depicted in Figure \ref{FigureInternals} (left panel).  The covariance matrix is $\hat{\Sigma}=\text{Diag}(\tilde{V})$, where $\tilde{V}$ denotes the vector of smoothed variances.
\item {\it Smoothed correlogram (SC).}  This option uses the same smoothed variances of SD but assumes that the correlation between $\varphi (A_{u})$ and $\varphi (A_{v})$ is stationary, {\em i.e.}~a function of the separation $|u-v|$.  Following a procedure in \citet*{Hall1994}, the sample correlogram is smoothed and then its negative Fourier frequencies are truncated to ensure positive-definiteness.  Figure \ref{FigureInternals} shows the sample and smoothed correlograms (right) for the FIAC data.
\item {\it Sampson-Guttorp Warping (SGW).}  Warping \citep{Sampson1992} does not assume stationarity in its smoothing of covariance.  Multi-dimensional scaling (MDS) is applied to the sample semivariogram, and then all but the first $k$ MDS components are dropped.  The top $k$ are re-expanded to produce a smoothed semivariogram, which is then translated to a covariance matrix.  In our experience, a minimum of $k=5$ components are necessary; we use $k=10$.
\item {\it Pseudo-inverse (PI).}  The sample covariance matrix cannot be directly used as $\hat{\Sigma}$ because it is singular when $F<U$; $\Sigma ^{-1}$ can be replaced with the pseudoinverse of $\Sigma$ \citep{Rao1962}, obtained by dropping the eigenvectors corresponding to eigenvalue $\lambda=$ 0 in the sample covariance matrix.
\end{enumerate}

\begin{figure}[!htbp]
  \includegraphics[width=1.55in]{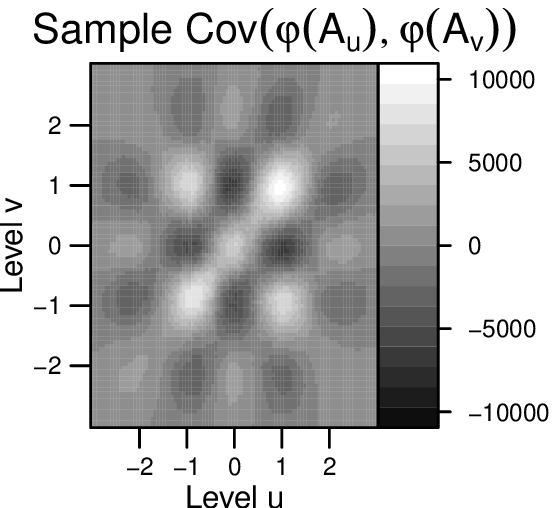}  \hspace*{0.2in}
  \includegraphics[width=1.3in]{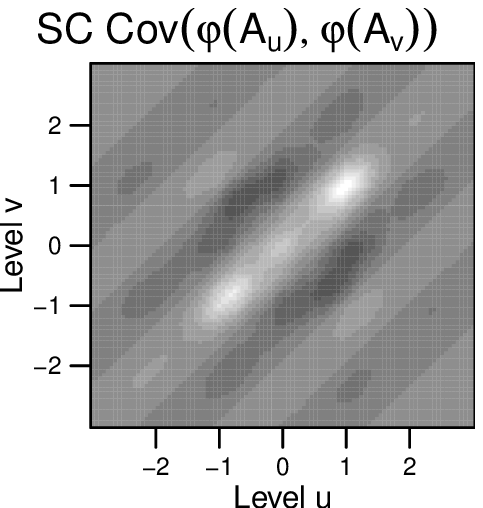}   \hspace*{0.2in}
    \includegraphics[width=1.3in]{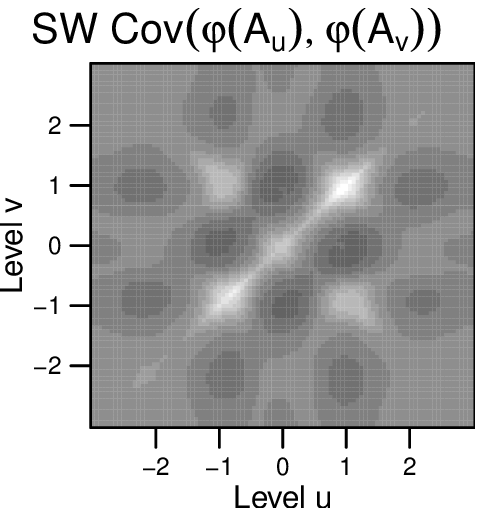}
  \caption{Covariance matrix of $\varphi(A_{u})$ estimated in three ways: the sample covariance (left), the SC method (middle), and the SGW method (right).}
\label{FigureCovariance}
\end{figure}

The first four options --- I, SD, SC and SGW --- represent a complexity continuum in modeling $C$.  Figure \ref{FigureCovariance} illustrates SC and SGW, comparing the results of SC (middle) and SGW (right) to the original sample covariance.  SC depends on the assumption that $\varphi (A_{u})$ is stationary, which accounts for the diagonal stripes in the SC covariance plot.  SGW is more flexible, yet we shall see that SGW tends to overfit $\Sigma$, which tends to increase the variance of the final regression estimates.

The fifth option, PI, unlike the others, does not estimate the function $C_\varphi(u,v)$ for arbitrary $u$ and $v$; it produces only $\hat{\Sigma}$, and only for a prescribed set of $u_{j}$.  While PI does incorporate both heteroskedasticity and correlation, we shall see that it too suffers from overfitting and yields highly variable regression estimates.

Our final recommendation is the SD covariance in the regression, as it provides the best balance between computation speed, flexibility and stability.  SC, SGW, and PI are flexible, however tend to overfit the sample covariance matrix, producing estimates of high variance.  We conduct a detailed comparison of the five approaches in Section \ref{SectionSimulationExperiment}.

\subsection{Design Selection}

\label{SectionDesignSelection}

The levels ${\bf u}=\{u_{j}\}$ are free to be chosen in LKC regression; both the number $U$ of levels and their locations (placement) need to be specified.  The primary tradeoff is between accuracy and speed.  If the levels are too few or too coarsely spaced, then estimation is fast but accuracy suffers.  On the other hand, if the levels are too many or too finely spaced, computational time is wasted with little marginal gain because the resulting ECs are highly correlated with one another.

A guiding principle is to minimize the theoretical variance of $\hat{\mathcal{L}}_{i}$, the regression-estimated LKCs.  Referring to (\ref{EquationModel}), let $X({\bf u})$ be the $U\times \dim (S)$ matrix of regressors, with $X_{ij}=\rho _{i}(u_{j})$, and denote the error covariance as $\Sigma ({\bf u})$ to emphasize its dependence on ${\bf u}$.  This gives
\begin{equation*}
V[\hat{\mathcal{L}}]=X^{\prime }({\bf u})\Sigma ({\bf u})^{-1}X({\bf u}).
\end{equation*}
Since $V[\hat{\mathcal{L}}]$ is a matrix, we require a real-valued function $g$ of $V[\hat{\mathcal{L}}]$ to perform optimization.  We consider A-optimality, where $g(M) = \text{tr}(M)$, and D-optimality, where $g(M) = \det(M)$ \citep{Pukelsheim2006}, and substitute an estimated covariance matrix from Section \ref{SectionErrorCovariance} for $\Sigma({\bf u})$.  Note that for a diagonal $\Sigma({\bf u})$, the optimizing $u_j$ lie at the extremes of the permissible range of ${\bf u}$, which requires $\Sigma({\bf u})$ to have nonzero off-diagonal terms for the optimization of ${\bf u}$ to be sensible.

\begin{figure}[!htbp]
  \includegraphics[width=2.2in]{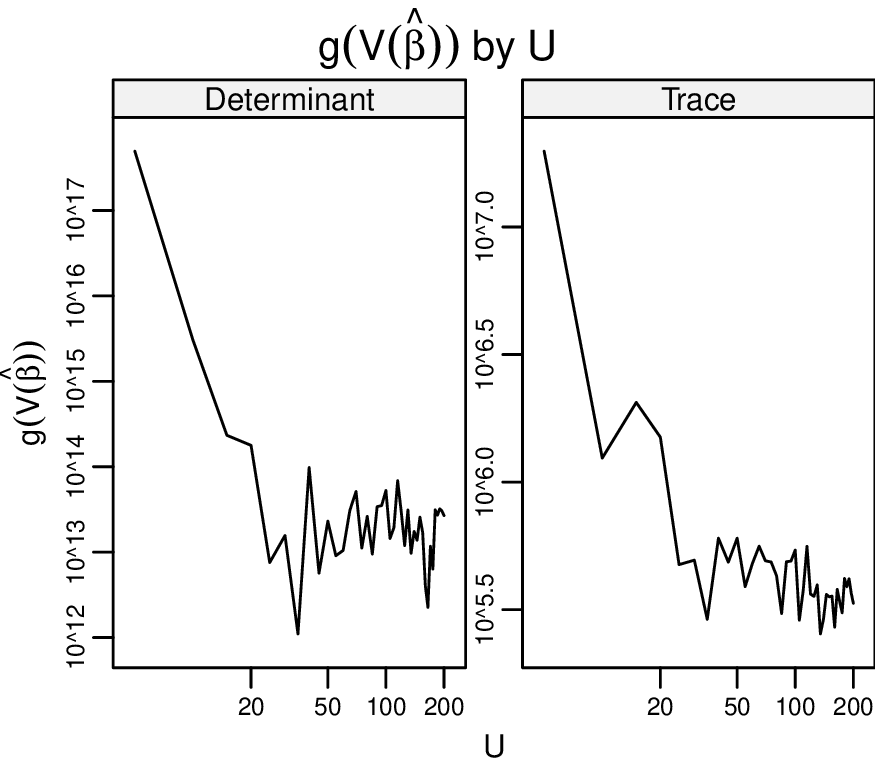}  \hspace*{0.4in}
    \includegraphics[width=1.5in]{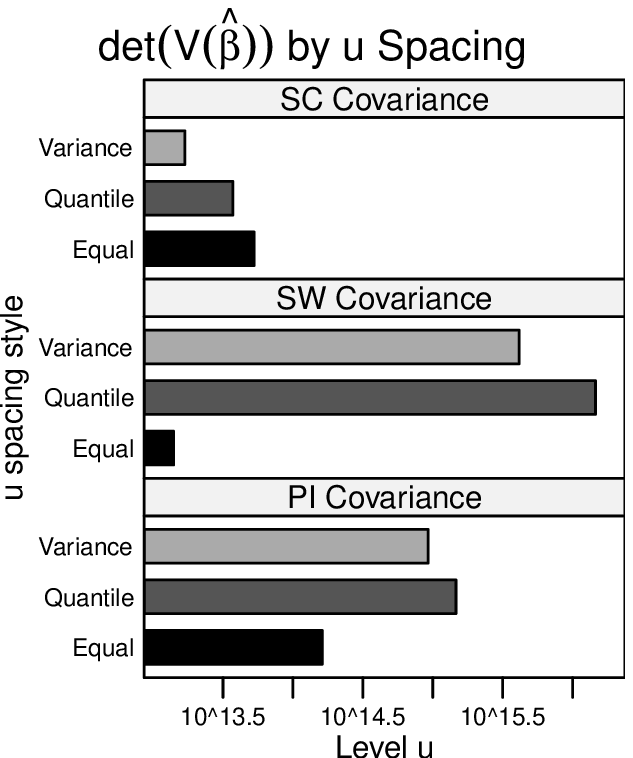}
  \caption{Factors of optimal design in the LKC regression: the number of levels $U$ (left) and the spacing mode of the levels ${u_j}$ (right).  The plots show how these factors impact the optimal design criterion $g(V[\hat{\beta}])$, for $g(M) = \det(M)$ and $g(M) = \text{tr}(M)$.  Level counts ranging from $U =$ 5 to 200 are considered along with equal, quantile and variance-based spacing.  The $U$ plot uses equal spacing and the spacing plot uses $U =$ 50 levels.}
  \label{FigureDesign}
\end{figure}

We first examine the question of how many levels to use.  Using the FIAC data, we tried possible values of $U$ between 5 and 200.  For illustration, we fix ${\bf u}$ to be $U$ equally-spaced levels between the field minimum and the maximum, though our results are not sensitive to this mode of spacing.  Figure \ref{FigureDesign} (left) shows the result for the SC covariance smoothing method as $U$ varies.  Although a larger $U$ may be better, the improvement levels off by about $U =$ 50.  Thus, we fix $U =$ 50 throughout our subsequent analysis.  The SGW and PI covariance options gave nearly identical results.

We next attempt to find an optimal distribution of levels.  Ideally, we could solve the optimization problem
\begin{equation*}
{\bf u}_{\text{optimal}}=\arg \min_{\bf u}g(X^{\prime }({\bf u})\Sigma ({\bf u})^{-1}X({\bf u})).
\end{equation*}
However, this search is challenging because the domain is $U$-dimensional.  In addition, the objective is discontinuous in $u$ due to a variety of instabilities in the covariance smoothing function.  Rather than attempting a direct optimization, we instead compare the objective values for a few placement heuristics:

\begin{enumerate}
\item {\em Equal spacing} of $u$ from $\min_{i,s}T^{(i)}(s)$ to $\max_{i,s}T^{(i)}(s)$.
\item {\em Quantile spacing} of $u$ at the 1\%, 3\%, $\ldots$, 99\% quantiles of the values in $T^{(i)}$, which places more levels where there are more field values.
\item {\em Variance spacing} of $u$, where placement density is proportional to $\hat{V}_u$, which places more levels where the EC is more variable.
\end{enumerate}

Figure \ref{FigureDesign} (right) shows the results for each of the options SC, SGW and PI.  Using SC, all design selection schemes perform similarly.  However, under the SGW and PI covariance matrices, equal spacing outperforms by a substantial margin.  We also experimented with the design styles in the simulation experiment of Section \ref{SectionSimulationExperiment} and found a similar advantage for equal spacing.  We henceforth use equally spaced levels exclusively.

\subsection{Comparison to Warping}

\label{SectionWarpingComparison}

The only published alternative method for LKC estimation is the warping approach of \citet{Taylor2007}.  The main idea of warping is to transform the input random fields such that they become roughly isotropic in a higher-dimensional space while their LKCs are held constant.

For comparison purposes, we provide a brief procedural review of the warping method.  Suppose the input random fields $T^{(i)}$ lie on a domain $S$ and have sample locations at $s_{k}\in S,\ k=1,\ \ldots,\ K$, and we are given a triangulation $\Delta$ of $S$ with a vertex at each $s_{k}$.  Warping transforms $\Delta$ to $\Delta^\prime \in \mathbb{R}^F$, a mesh in $F$-dimensional space also with $K$ vertices.  The coordinates of each vertex of $\Delta^\prime$ are the observed field values $T^{(1)}(s_k),\ T^{(2)}(s_k),\ \ldots\ T^{(F)}(s_k)$ from a particular location $s_k$.  Formally,
\begin{equation*}
  \text{$j^{\text{th}}$ coordinate of $k^{\text{th}}$ vertex of $\Delta^{\prime }$} = \frac{T^{(j)}(s_{k})}{\sqrt{\sum_{i=1}^{F}T^{(i)}(s_{k})^{2}}}
\end{equation*}
The (Euclidean) intrinsic volumes of $\Delta^{\prime }$ are then the estimates of the LKCs $\mathcal{L}_{1},\ \ldots ,\ \mathcal{L}_{\dim (S)}$.  The precise formulation is provided in \citet{Taylor2007}, and involves the edge lengths, interior angles, surface areas and volumes of the triangulation in the transformed space.

An apparent downside of warping is its conceptual as well as computational complexity.  The geometric calculations in the transformed space can be time-consuming for a large triangulation.  Unlike $\Delta$, $\Delta^{\prime }$ does not lie along a regular grid, so there are no cost-saving computational measures available and runtimes can be up to an order of magnitude longer, which makes the method time-prohibitive when multiple sets of LKCs need to be found ({\em e.g.}~for a time series of repeated experiments).  In contrast, our LKC regression requires the computation of only the EC instead of all the intrinsic volumes.  It also works directly in the space of the input fields, where the regular layout of sample locations permits speedy computation.

LKC regression is fundamentally more transparent and conceptually accessible, while the warping procedure may be more challenging for nonspecialists to understand.  Diagnostics are difficult under warping because it can be unclear to detect the cause of an erroneous estimate, while LKC regression is interpreted simply as the curve that best matches expected with observed ECs.  Any necessary diagnostics --- such as outliers or data integrity issues --- are revealed in a straightforward manner by examining the EC profile.

\section{Simulation Study}

\label{SectionSimulationExperiment}

To test the accuracy of the Lipschitz-Killing curvature (LKC) estimates resulting from our regression procedure proposed in Section \ref{SectionLKCRegression}, a numerical experiment using simulated random fields was performed.  Our aims are twofold: first, to determine optimal settings for our LKC regression; and second, to compare it to the main alternative, the warping method \citep{Taylor2007}.

Rather than on the LKC estimates themselves, our simulation study focuses on the estimated 95\% threshold for the field maximum, which are determined by the LKC estimates via (\ref{EquationThreshold}).  This threshold is the driver for 95\% hypothesis tests, and with practitioners' needs in mind, the primary purpose and application for which our LKC method was developed.  Our simulation study computes the standard deviation and bias in this threshold.  To assess the generality of the method, we test a range of random fields commonly found in practice.

\subsection{Experimental Metrics}

\label{SectionExperimentalMetrics}

We compare the LKC estimation methods by their computational runtime, error and bias.  
Care has been taken in all cases to employ the optimal implementation of each method and to exclude fixed costs ({\em e.g.}~loading the data).  Standard deviation describes how much an estimated threshold $\widehat{u_{95\%}}$ varies over multiple simulations.

Bias reveals how far, on average, $\widehat{u_{95\%}}$ is from the true threshold.  Bias is harder to measure as it requires a specification for the input random fields --- a domain and covariance function --- for which the true threshold is known.  For stationary and isotropic fields, it is known that the LKCs can be comparatively simple to compute, particularly if $S$ has a simple topology.  We thus base our simulations on the Gaussian covariance ($C(x,y)=e^{-\alpha |x-y|^{2}}$) on the square, cube and sphere.
The true LKCs and thresholds, calculated via the ECH, are shown in Table \ref{TableTruth} (first row).  Our experiment is made up of simulations from these three field types
with suitable values of $\alpha$.  The first LKC, $\mathcal L_0$,  is not included, Since $\mathcal L_0$ is the Euler characteristic, it is known deterministically, and, for the three examples in the table it is 1, 1, and 2, respectively.

\begin{table}[!htbp]
  \begin{center}
    \begin{tabular}{ccc|cccc|ccc}
      \multicolumn{3}{c|}{\textbf{Square}} &
      \multicolumn{4}{|c|}{\textbf{Cube}} &
      \multicolumn{3}{|c}{\textbf{Sphere}} \\
      $\mathcal{L}_1$ & $\mathcal{L}_2$ & $u^{\ast}_{95\%}$ &
       $\mathcal{L}_1$ & $\mathcal{L}_2$ & $\mathcal{L}_3$ &
      $u^{\ast}_{95\%}$ &  $\mathcal{L}_1$ & $\mathcal{L}_2$ &
      $u^{\ast}_{95\%}$ \\ \hline
       $2 \sqrt{2\alpha}$ & $2\alpha$ & $-$ & $3 \sqrt{2\alpha}$ &
      $3 (2\alpha)$ & $(2\alpha)^{3/2}$ & $-$ &  $0$ & $4\pi (2\alpha)$ &
      $-$\\
       28.3 & 200 & 3.72 & 19.0 & 120 & 253 & 3.96 &  0 & 503 & 3.96
    \end{tabular}
  \end{center}
  \caption{True continuous LKCs and 95\% thresholds for continuous random fields in the square $(\alpha = 100)$, in the cube $(\alpha = 20)$, and on the surface of the sphere $(\alpha = 20)$.  Both analytical expressions (upper row) and simulation-specific numerical values (lower row) are provided.  In each case, $L_0$ has an easy known form as the EC of the domain $S$.  The other LKC formulae are based on analytical results from \citet{RFG}.  The 95\% thresholds do not have easy analytical forms, but are calculated using the ECH.  These expressions act as the ``truth'' in our experimental assessment of bias.}
  \label{TableTruth}
\end{table}

A complication for bias computation is that the known thresholds are for continuous random fields, whereas our experimental fields are generated at discrete sample sites (see Section \ref{SectionObservedExcursionSets}).  Discrete and continuous random fields have different thresholds even when they share the same covariance function and domain.  There are no known analytical forms for the LKCs of discrete fields \citep{RFG}, but it is not hard to show that, for the random fields under consideration, the threshold for discrete fields converges to that for the continuous limit as the simulation resolution increases.  The convergence of discrete to continuous thresholds provides a basis for computing bias; for each estimation method, the average $\widehat{u_{95\%}}$ at varying grid sizes thus forms a convergent pattern.  The approximate limit yields an estimate of the continuous threshold, which can be compared to the known true value to measure bias.  (Figure \ref{FigureThresholdMedian} provides an illustration.)

\subsection{Experimental Factors}

\label{SectionExperimentalFactors}

Observed random fields are commonly classified by the first four factors listed in Table \ref{TableExperimentalFactors}: number of fields, domain, covariance function, and resolution.  Our experiment is factorial and tests 10,000 replicates of each combination of these factors.  The domains are the unit square, the unit cube, and the surface of a unit sphere.  The grid sizes range from $G$=5 to 200 sample locations in each direction, which spans most applications and reveals trends in the estimates as the resolution increases.  We simulate $F$=15 Gaussian random fields in each run, representing a typical thresholding experiment.  Each field is equipped with Gaussian covariance for a suitable value of $\alpha$.  (All fields are simulated via the turning bands method \citep{Mantoglou1982}: we use 1,000 lines and 4$G$ points along each line, twice as many as the authors recommend to help ensure that simulation error is small compared to estimation error.  We adapt the implementation of \citet{Schlather2009}.)

To determine the best settings for LKC regression, we also test the five positive-definite covariance construction options discussed in Section \ref{SectionErrorCovariance}.  All estimations use $U=50$ equally-spaced levels $u_{j}$.  For the spherical example $\mathcal{L}_{1}$ is always zero for topological reasons, and so neither  the LKC regression nor warping  attempt to estimate it.  We employ the technique of \citet{Taylor2007} for warping (for the triangulation $\Delta$, we input a mesh of standard right triangles for the cube and plane and latitudinal-longitudinal triangles for the sphere).

\begin{table}[!htbp]
  \begin{tabular}{ll}
    \textbf{Factor} & \textbf{Values} \\ 
    \# fields & $F = 15$ \\ 
    Domain & Square, cube, sphere \\ 
    Covariance & $c(x, y) = e^{-100|x-y|^2}\ \text{(square)},\ e^{-20|x-y|^2}\ 
    \text{(cube and sphere)}$ \\ 
    Grid size & $G = 5,\ 10,\ 20,\ 50,\ 100,\ 200\ \text{(square only)}$ \\ 
    Method & Lipschitz-Killing curvature regression (LKCR), warping \\ 
    $\Sigma$ (LKCR only) & From Section \ref{SectionErrorCovariance}: I, SD, SC,
    SGW, PI \\ 
    \# levels (LKCR only) & $U = 50$ \\ 
    & 
  \end{tabular}
  \caption{Factors for the simulation experiment of Section \ref{SectionSimulationExperiment}.  Factors describe both the nature of the input random fields (first four rows) and the LKC estimation methods (last three rows).}
  \label{TableExperimentalFactors}
\end{table}

\subsection{Results}

\label{SectionResults}

Overall, we find that LKC regression is much faster than warping, and maintains similar accuracy.  The increase in speed is roughly of a factor of up to eight for high-resolution random fields on the square and cube, which are the most common domains in applications.  This gain comes with negligible loss in standard deviation or bias.  Among LKC regression covariance styles, the smoothed diagonal (SD) method is a consistent winner, providing the best combination of accuracy and speed.

\subsubsection{Runtime}

\label{SectionRuntime}

\begin{figure}[!htbp]
\includegraphics[width=4.5in]{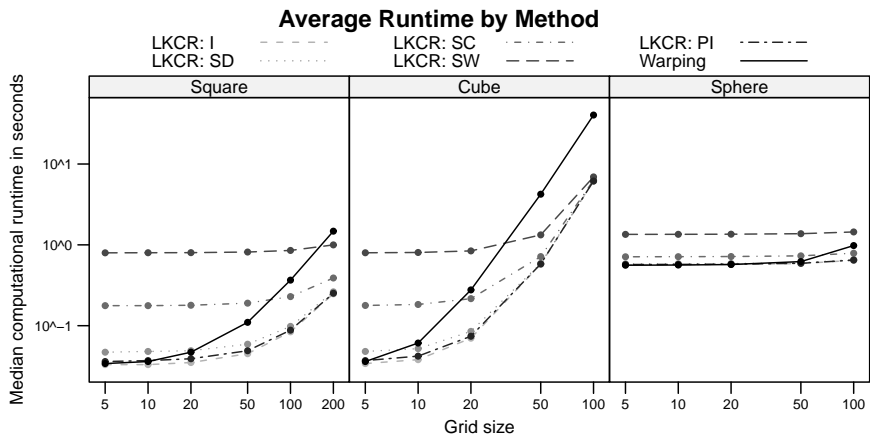}
\caption{Median runtimes for estimating LKCs in the simulation experiment.  Results are given by method (line styles), random field domain (square, cube or sphere) and simulation grid size (horizontal axis).  Medians are over 10,000 independent replications of the input data.}
\label{FigureRuntime}
\end{figure}

Computational runtime provides a compelling case for LKC regression.  The averages shown in Figure \ref{FigureRuntime} show a dramatic advantage for most LKC regression methods (dashed and dotted lines) against warping (solid line).  This advantage is most pronounced for large grids that are commonly found in practice.  For cubic fields with $G =$ 100, warping takes an average 40 sec compared to 6 sec for LKC regression methods I, SD and pseudoinverse (PI).

Among LKC regression styles, the slowest are smoothed covariance (SC) and Sampson-Guttorp warping (SGW).  These two options involve complex estimation of the off-diagonal terms in the regression covariance matrix, which thus incurs a heavy fixed cost in covariance estimation that dominates the runtime.

\subsubsection{Standard Deviation}

\label{SectionStandardError}

\begin{figure}[!htbp]
  \includegraphics[width=4.5in]{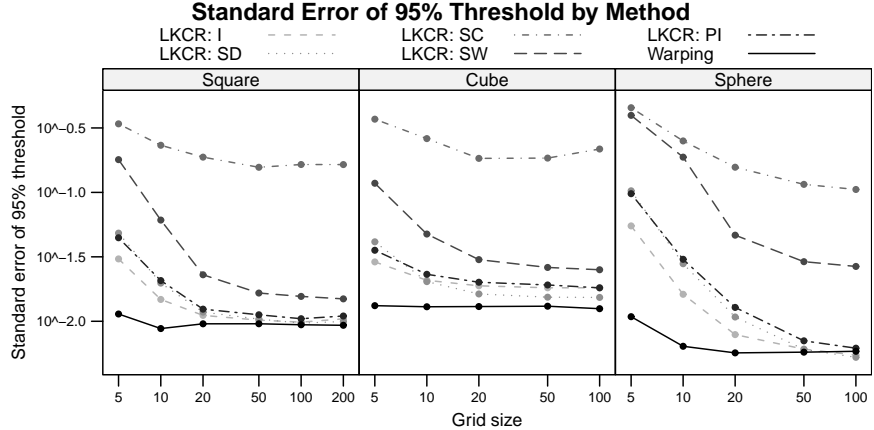}
  \caption{Standard deviations of 95\% threshold estimates.}
  \label{FigureThresholdSE}
\end{figure}

Comparing the standard deviation of LKC regression and warping method depends on simulation resolution.  Figure \ref{FigureThresholdSE} shows generally similar standard deviations between warping and the top-performing LKC regression styles, identity and SD covariances, for large grid sizes.  For fields with $G =$ 50, the differences are
slight: 0.0096 for warping versus 0.0103 for LKC regression using SD covariance for square fields, 0.0056 versus 0.0059 for spherical fields, and 0.0127 versus 0.0158 for cubic fields.  Moreover, both standard deviations are negligible when viewed with regard to the underlying application.  Applied random fields usually have thresholds
of around 5 (see Section \ref{SectionApplication} for an example).  Standard deviations on the scale of those we observed have little impact on the threshold.  They essentially do not affect the task of identifying significant regions in high-resolution fields.  Warping has smaller standard deviations for fields with low resolutions, particularly $G \leq$ 10.

The best LKC regression covariance structure is shown to be SD, although I and PI also perform well.  The poor performance of SC and SGW is likely due to overfitting; both fit nonparametric models of the regression covariance matrix.  Errors in these procedures accumulate in the ultimate regression, causing higher standard deviations.  The I and PI covariance models are less flexible but more parsimonious, stabilizing the estimates.  SD is similar to I, but comprises a relatively conservative smoothing procedure (see Figure \ref{FigureInternals}) to model heteroskedasticity.

\subsubsection{Bias}

\label{SectionBias}

\begin{figure}[!htbp]
  \includegraphics[width=4.5in]{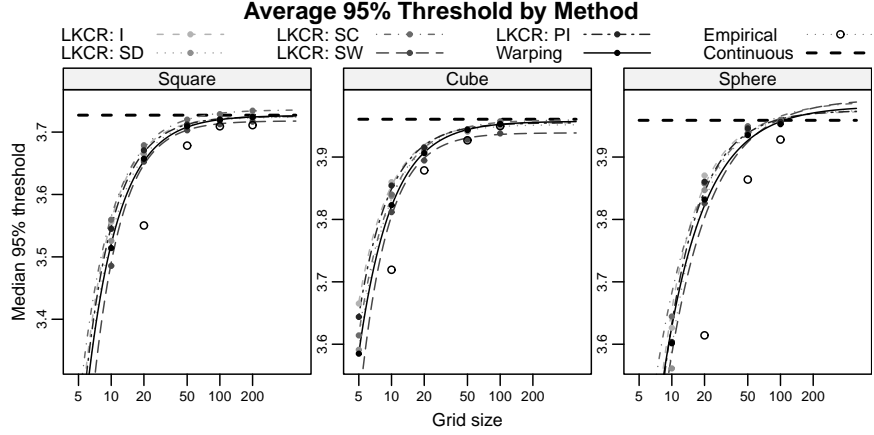}
  \caption{Medians of 95\% threshold estimates.  The extrapolated curves are fitted using nonlinear regression, Equation (\ref{EquationLimit}).  The unlinked empirical dots are the ``true'' thresholds approximated using 150,000 simulated discrete random fields.}
\label{FigureThresholdMedian}
\end{figure}

As mentioned in Section \ref{SectionExperimentalMetrics}, a direct bias calculation is difficult because true values of LKCs for discrete random fields are unknown; true values are known only for continuous random fields with Gaussian covariance.  However, the continuous LKCs can be approximated by extrapolating from a pattern of discrete LKCs at varying resolutions.  Similarly, the median estimated 95\% thresholds in Figure \ref{FigureThresholdMedian} show clear convergent patterns.  Assuming convergence of order $\varsigma$, the limit is found using nonlinear regression: 
\begin{equation}
u_{95\%}(g)=u_{95\%}^{\ast }+\beta g^{\varsigma }+\varepsilon (g).
\label{EquationLimit}
\end{equation}
Here, $u_{95\%}^{\ast }$ denotes the implied continuous threshold, while $u_{95\%}(g)$ is the estimated 95\% threshold at grid size $g$.  The unknown parameters $\beta$ and $\varsigma$ govern the rate of approach to the continuous value, for which typical estimates are $\hat{\beta}=-7.0$ and $\hat{\varsigma}=-1.5$; $\varepsilon$ is an error term that we take as normal for simplicity.

The end result of regression (\ref{EquationLimit}) is $\widehat{u^{\ast}_{95\%}}$, an approximate continuous threshold.  One regression is needed for each method and random field type, representing a single pattern of dots in Figure \ref{FigureThresholdMedian}.  The fitted regression lines are also shown; bias is measured as the difference between approximate and true continuous thresholds.

\begin{figure}[!htbp]
  \includegraphics[width=4.5in]{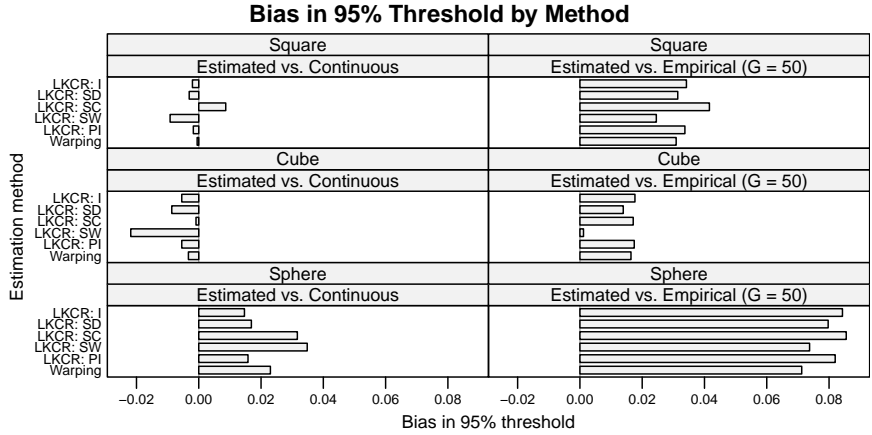}
  \caption{Biases of two types of 95\% threshold estimates: the continuous threshold estimated via extrapolation (left); and the discrete threshold for resolution $G =$ 50.}
  \label{FigureThresholdBias}
\end{figure}

Neither LKC regression nor warping emerges as the clear winner in bias.  Figure \ref{FigureThresholdBias} (left column) indicate small, similar levels of bias across all methods.  The only outlier is the SGW covariance, which again produces poor estimates, likely to be the result of an overfitting of covariance.

An alternative approach to measuring bias is to approximate the true 95\% thresholds for discrete random fields with intensive simulation.  Given $B$ simulated fields, the 95\% sample quantile of the field maxima provides an approximation (we use $B=$150,000).  These estimates appear in Figure \ref{FigureThresholdMedian} (\textquotedblleft Empirical\textquotedblright), and are used to measure bias in Figure \ref{FigureThresholdBias} (right column).  The result is the same: all methods demonstrate a similar level of bias, however Figure \ref{FigureThresholdMedian} also reveals an interesting self-regulatory mechanism in LKC estimation.  The estimated discrete thresholds always lie between the continuous thresholds computed via the ECH (\ref{EquationThreshold}) and the discrete thresholds computed via simulation.  Despite their basis in continuous theory, the estimated thresholds of the LKC regression approach the discrete truth.

\section{Applications}

\label{SectionApplication}

\begin{figure}[!htbp]
  \includegraphics[width=1.4in]{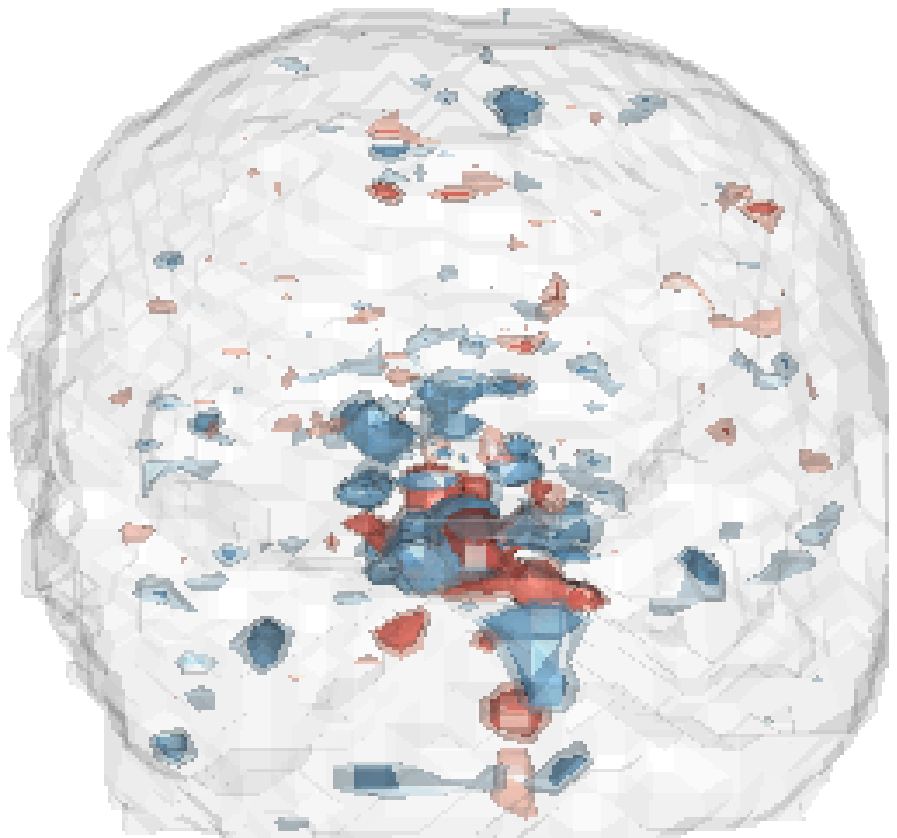}  \hspace*{0.2in}
    \includegraphics[width=1.4in]{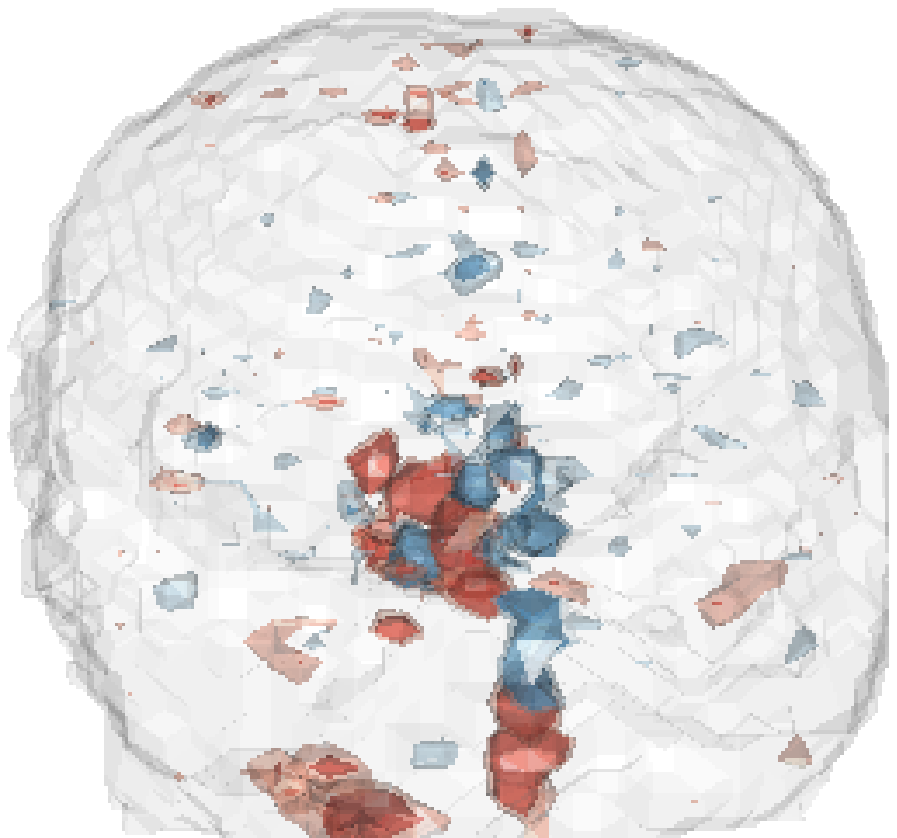}  \hspace*{0.2in}
      \includegraphics[width=1.4in]{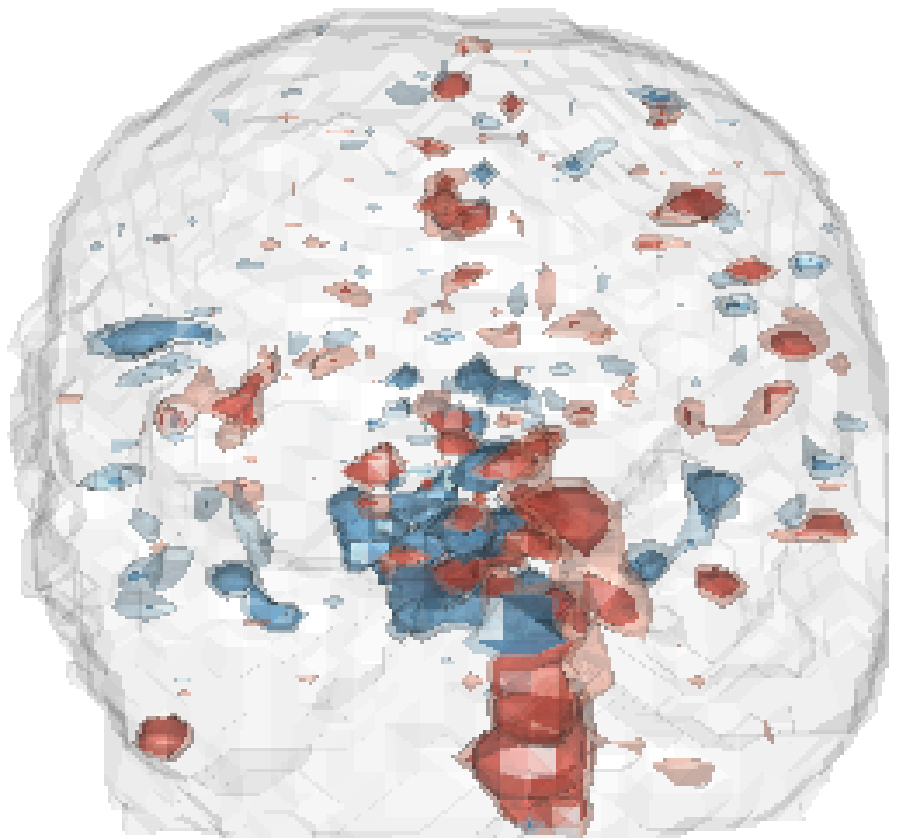}
  \caption{Examples of three residual fields from a voxel-wise AR(1) model fit to the FIAC brain data.  The scans shown come from the first experimental subject at times 1, 50 and 100.  The outer gray layer shows the surface of the brain, while the red and blue show domains with positive and negative residuals, respectively.  Light and dark colors represent the data thresholded at the single t-test 95\% and 99\% thresholds of 1.96 and 2.58.}
  \label{FigureApplicationExample}
\end{figure}

In this section we discuss two applications. One, on the analysis of CMB data, is discussed  briefly at the end of the section. For a more detailed example, we return to the fMRI data discussed in the Introduction.

\subsection{A fMRI example}
We return to the language priming experiment of \citet{Dehaene2006} and apply our LKC regression method.  

Contrary to the group-level comparison of Section \ref{SectionIntroduction}, we conduct a within-subject study to identify stimulated regions for each specific individual.  Here the stimulus refers to the difference in activation between the ``different speaker'' and ``same speaker'' conditions, under which each subject hears a sentence repeated twice either by two different speakers or by the same speaker.  The subject is then scanned 191 times, alternatively undergoing the two stimuli.  There are 14 such data sets, one for each individual, which are analyzed by \citet{Taylor2006Inference} using the warping thresholding method.  We reanalyze the data using our LKC regression method.

Prior to the analysis we conducted a standard preprocessing \citep[as detailed in Section 2 of][]{Taylor2006Inference} of the data: First, a subject's scans are aligned to correct for motion drift \citep{Smith2004}.  An AR(1) time series model is then fitted to every voxel of the brain.  The response variable is the fMRI activation level, which is a time series indexed by each of the 191 time points; in the AR(1) model, the stimulus effect is treated as a regression coefficient in front of a time-varying contrast that represents the level of stimulus administered at each time.  There is one such model for every voxel in the brain; fitting all such models produces a random field of estimated stimulus effects across the brain: this is the field of interest, analogous to the field $T$ in Section \ref{SectionIntroduction}.  The random fields $T^{(i)}$, $i = 1,\ \ldots,\ 191$ correspond to the AR(1) residuals at every voxel and time.  As in Section \ref{SectionIntroduction}, the threshold is obtained from the $T^{(i)}$ and then applied to $T$ to identify regions where the stimulus effect is significant.  There is one analysis --- one $T$ and one set of $T^{(i)}$ --- for every subject.  As an example, three residual fields for the first subject in the study, corresponding to three separate time points, are shown in Figure \ref{FigureApplicationExample}.  In \citet{Taylor2006Inference}, warping is used to estimate the threshold, which found values of around 5 for each subject.

We are interested in analyzing how these thresholds, found using warping \citep{Taylor2007}, compare to those found using our proposed LKC regression method.  For the purpose of such a comparison, the thresholds are estimated using warping and our LKC regression with smoothed diagonal (SD) covariance (which is found to be the best-performing option found in Section \ref{SectionSimulationExperiment}).

\begin{table}[!htbp]
  \begin{center}
    \begin{tabular}{r|rr|rr|rr|rr|rr}
      & \multicolumn{2}{c}{\textbf{95\% thres.}} &
      \multicolumn{2}{|c}{$\boldsymbol{\mathcal{L}_1}$} &
      \multicolumn{2}{|c}{$\boldsymbol{\mathcal{L}_2}$} &
      \multicolumn{2}{|c}{$\boldsymbol{\mathcal{L}_3}$} &
      \multicolumn{2}{|c}{\textbf{Runtime}}\\
      \textbf{} & LKC & Warp & LKC & Warp & LKC & Warp
      & LKC & Warp & LKC & Warp\\ \hline
            1 & 5.01  & 4.94  &1,686  & 274  &7,774  &2,343  &26,706  &29,382  &18.3   &38.7 \\      
      2 &  4.93 & 4.75 & 1,237 & -32  & 7,801 & 2,351  &23,878  & 32,334  & 12.5  & 34.9 \\           
      3 &  5.10 & 4.98 & 943  &  201 &4,694  & 2,242 & 15,234 & 30,659  & 15.0  & 38.7 \\           
      4 &  5.03 & 5.14   & 1,252  & 249 & 5,584 & 2,169  & 19,886  & 29,537  & 15.8  &33.2 \\           
      5 & 5.06  & 5.11  & 1,961  & 49 & 9,484 & 2,277  & 30,229 & 30,683 & 12.7  & 33.8 \\      
      6 &   4.99 &4.98  &1,158  & 22 & 7,712  & 2,427 & 23,317 & 34,124 & 14.8  &32.7 \\           
      7 &  5.03  & 5.09   &1,263  & 419  & 4,702 & 2,091  &  17,543  & 28,103  & 16.7  & 37.5 \\           
      8 &  4.94  & 5.05  & 1,737  &  382  & 8,860  & 2,162  & 28,163 & 26,701  & 13.9  & 34.1\\      
      9 &  4.97 & 5.02  & 2,017  & 191 & 8,782  & 2,235 & 31,009  &31,327  & 20.9   & 35.4 \\      
      10 &  5.15 & 5.28   & 1,437  &   386   & 4,857   &2,245  &  17,988  & 30,403  &18.5   &  42.6 \\           
      11 & 4.98  & 4.91 & 2,394  & 320  & 8,967  & 2,274  & 30,473 & 31,133 & 15.3  & 37.3 \\           
      12 &  5.01  & 4.89  & 1,333  & 307 & 6,528 & 2,113  &21,540  & 28,898 & 15.3  & 36.1 \\      
      13 & 5.01   & 5.01  & 1,933  & 255 & 7,658  & 2,063  & 27,509 & 29,670  & 13.6  & 34.5\\      
      14 &  5.01 & 5.13  & 1,757  & 105  & 7,724 & 2,262 & 26,367 & 32,376  & 12.5  & 34.1 \\           
    \end{tabular}
  \end{center}
  \caption{Per-subject LKC estimation results for LKC regression (LKCR) and warping.  The data are residual random fields from a voxel-wise AR(1) model fit to the FIAC brain scan data, with $F =$ 191 fields per subject.  The estimated 95\% thresholds (left), LKCs (middle) and runtime (in seconds, right) are displayed for each subject and method.}
  \label{TableApplication}
\end{table}

Table  \ref{TableApplication} compares the LKC estimates and the computational time of each method.  The major difference is in the latter (rightmost column): the LKC regression is over twice as fast as warping. At first, while any improvement in computational speed is desirable, it would not seem that this improvement is all that important. Statistical analysis of fMRI data is typically done offline, as opposed to the actual production of fMRI data which is online and must be minimised for subject comfort and cost considerations. 

Nevertheless, computation time  is still an issue, as evidence by the fact that
 practitioners have developed data structures optimized for 3D fields \citep[e.g.,][]{Theis2005} and schemes for parallel processing \citep[e.g.,][]{Zhao2007,Wilde2009}.  Computational firepower is necessary because over the course of an fMRI study, the number of hypothesis tests conducted and thresholds found easily reaches into the hundreds.  These computational demands occur for two reasons: first, there are usually multiple experimental stimuli, each with a categorical set of possible conditions (for instance, \citet{Taylor2006Inference} consider three factors with two or three conditions each, yielding 12 significance analyses); and second, it is also common to perform subject-specific tests (one test per subject), or time-specific tests (one test per unit time) \citep[see][]{Beckmann2003}.  In light of these repeated computational demands that are inherent to various applications, our LKC regression can save substantial development time.

In addition to speed, accuracy also carries great importance for practictioners.  As seen from the spread of the estimated 95\% thresholds (leftmost column), LKC regression and warping produce very similar thresholds, with the mean values of the 14 subjects being 5.01 and 5.02, respectively.   Thus, as to be expected, the activated regions change only slightly (not shown) and we are still able to reproduce the  activation in the left and right mid-temporal gyri seen in \citet{Taylor2007}. 

What we have not yet checked in detail, but, from preliminary studies are quite certain is true, is that in higher dimensional problems the gain in speed of LKC regression over warping will be quite significant. What is definitely true is that the LKC approach, even in this three dimensional example, involves much less coding, since, following Section 3, it is basically just a regression analysis.

As a side note, it is interesting that the similarity in thresholds between LKC regression and warping belies often large differences in the estimated LKCs (middle three columns).   These phenomena underscore the fact that different LKC combinations can lead to the same 95\% thresholds, which are what matter in practice. It is impossible to tell which is correct because only the thresholds are observable, not the LKCs themselves. Interestingly, the discrepancies in $\mathcal{L}_{1}$ and $\mathcal{L}_{2}$ are not apparent in the cubic random field simulations of Section \ref{SectionSimulationExperiment}, where LKC regression and warping produce approximately the same LKC estimates. We suspect that isotropy plays a role: the simulated cubic fields are isotropic, while the FIAC data are likely anisotropic.

\subsection{A cosmic microwave background radiation example}
Working from an earlier version of the current paper, 
\citet{Marinucci-LKC} have used LKC regression to estimate the LKCs (`Minkowski functionals' in their language) for   CMB fluctuation  models. 

To quote from their paper: ``A general trend in modern cosmological research is the implementation of more and more sophisticated statistical
tools to perform data analysis. Indeed, as well-known cosmological data have reached over the last decade an un-
precedented accuracy, so that it has become customary to speak about a golden era for Cosmology, featuring a data
deluge from a bunch of satellite - and ground based-experiments. As the data grow in size and precision, more and
more detailed questions can be addressed, and exploiting techniques at the frontier of statistical and mathematical
research becomes mandatory to warrant a full exploration of the available evidence. 

Among these techniques, stochastic geometry tools have now become very well established, especially in the field of
Cosmic Microwave Background radiation experiments. In this area, one of the most popular geometric tools for data
analysis are certainly the so-called Minkowski functionals (MFs), which have been extensively exploited as tools to
search for non-Gaussianities, anisotropies, asymmetries and other features of CMB data. The use of MFs in Cosmology
goes back at least to [1, 2]; a complete bibliography would certainly include hundreds of entries, so we refer only to
the earlier works by [3--10] and to the more recent ones by [11--16]." (References not included here.)

We recommend their paper, which contains detailed calculations and scientific discussions, as an excellent application of the techniques developed here (and, of course, for the many other things there).

\section{Conclusion}

\label{SectionConclusion}

We have presented a new method, the Lipschitz-Killing curvature regression, for estimating tail probabilities of a Gaussian random field.  The LKC regression procedure is to estimate the Lipschitz-Killing curvatures and then substitute them into the Euler characteristic heuristic (\ref{EquationThreshold}) to generate a $(1-\alpha)$ threshold.  The estimation does not require knowledge of the covariance structure of the field or even whether it is isotropic.  It need not be Gaussian, although it must be Gaussian related. The high efficiency of our method hinges on the fact that the Euler characteristic can be well estimated from the data when the exceedence level $u$ is small or moderate.  The LKC regression allows us to leverage the estimation
strength at a low level $u$ for accurate approximation at high exceedence level.

The primary advantage of the LKC regression is its simplicity and speed.  It is transparent, which offers straightforward interpretation, implementation, and diagnosis.  It runs faster than its chief competitor, warping \citep{Taylor2007}, with comparable accuracy.  The key procedural difference is that LKC regression computes Euler characteristics of the input random fields, which is a fast computation when the sample sites lie on a grid.  In practical terms, for the two- and three-dimensional random fields we considered, the gain in speed ranges from a factor of two to a factor of eight.  It is most pronounced for high resolution fields, which are commonly found in practice.  We believe that the gain rate will be considerably higher in higher dimensional cases, such as the scale-space (5D) of \cite{S-WOR95} and the rotation-space fields (8D) of \cite{Shafie:Sigal:Siegmund:Worsley:2001}, for both of which the analytic evaluation of the LKCs is orders of magnitude harder than the isotropic cases we have considered in this paper (cf.\ \cite{AST}).

\section*{Appendix}
\label{SectionAppendix}

\setcounter{equation}{0} \renewcommand{\thesection}{A}

In this section we shall briefly discuss the Gaussian kinematic formula, which gives an
exact formula for the
expected Euler characteristic $E[\varphi(A_u)]$. We shall explain its
general structure without going into technical details.

\subsection{Gaussian and Gaussian related random fields}

The basic building blocks of all the random fields we consider are
smooth Gaussian random fields, or processes,
\beqq
g \ = \ \left(g^1,\dots,g^k\right)\ \:S\in \real^N\to\real^k.
\eeqq

The first assumption that we place on $g$, and the most important one,
is that each of its $k$ components are twice differentiable, and that
these derivatives are themselves continuous. Some additional minor
assumptions of non-degeneracy also need to be made, but since these
almost always hold in practice we direct the interested reader to
Chapter 11 of \cite{RFG} or Chapter 4 of \citet*{ARF} for details. For
the rest of this paper we shall assume these conditions are met.  The second
assumption is a minor one, that all means be fixed at zero. A little
more restrictively, we also assume that the $g^j$ have constant
variance throughout $S$. Note that this is a much weaker assumption
than either isotropy or stationarity, which we do not require, and is
achievable in general by replacing a random field that does not have
constant variance by a normalized version of itself. Unless stated
otherwise, for the remainder of this section we shall assume that this
constant variance be 1. Finally, we assume that the components of $g$
are all independent and identically distributed.

Real valued Gaussian related random fields $f$ are defined by taking a
smooth function $F\:\real^k\to\real$ and setting $f(s) = F(g(s))$ for
all $s\in S$, where $g$ is as above. Gaussian related fields are
typically quite different to Gaussian ones (for which $k=1$ and $F$ is
the identity function). Three useful examples are given by the
following choices for $F$, where in the third we set $k=n+m$.
\beq
\label{rob-chi}
\sum_1^k x_i^2,\qquad
\frac{x_1\sqrt{k-1}}{(\sum_2^k  x_i^2)^{1/2}},\qquad
 \frac{m\sum_1^n  x_i^2}{n\sum_{n+1}^{n+m}  x_i^2}.
\eeq
The corresponding random fields are known as $\chi^2$ fields with $k$
degrees of freedom, the $t$ field with $k-1$ degrees of freedom, and
the $F$ field with $n$ and $m$ degrees of freedom. These three random
fields all have very different spatial behavior, and each is as
fundamental to the statistical applications of random field theory as
is its corresponding univariate distribution to standard statistical
theory.

Throughout we shall assume that $F$ is twice continuously
differentiable, and that for all real $u$ the sets
$F^{-1}[u,\infty)\subset \real^k$ are well behaved, as
  described in the following subsection. Under these conditions,
  Gaussian related fields are also twice continuously differentiable,
  with well behaved excursion sets.

\subsection{Parameter spaces}

The modern general theory of Gaussian fields, as developed in
\cite{RFG}, allows the parameter $S$ to be what is known as a Whitney
stratified manifold, satisfying some mild side conditions.  Roughly
speaking, these are compact subsets of $\RN$ which can be written as a
disjoint finite  union
\beqq
S= \bigsqcup_{i=1}^N \partial_iS,
\eeqq
where $\partial_iS$ is an $i$-dimensional manifold (the open
$i$-dimensional `boundary' of $S$).

An easy example is given by a $N$-dimensional rectangle, in which
$\partial_N S$ is its interior, $\partial_{N-1} S$ the collection of
its (open) $(N-1)$-dimensional faces, etc., down to $\partial_0 S$
which is the collection of its corners. If $S$ were a ball, we would
decompose it into its interior, again $\partial_N S$, and its surface,
the sphere $\partial_{N-1} S$.

The `Whitney' part of the definition contains rules about how all the
various pieces must be put together, and details can be found  in
\cite{RFG}.  More or less every (non-fractal) parameter space arising in 
statistical practice  satisfies these rules.

\subsection{The expected Euler characteristic of excursion sets}

The Euler characteristic $\p (S)$ of a well behaved set $S\subset \RN$ is 
a topological invariant, that, typically, is easy to
compute. For example, if $S$ is two-dimensional, then $\p (S)$ is
simply the number of connected components in $S$ minus the number of
holes. If $S$ is three-dimensional, then it is the number of
components minus the number of handles plus the number of internal
holes.

If $S$ is made up of a union of cubes, all of whose corners sit at the
points of a rectangular lattice, then
\beqq
\p (S)\ = \ \sum_{n=1} (-1)^n \mu_n, 
\eeqq
where $\mu_n$ is the number of {\it distinct} facets of dimension $n$
in $S$.

The central result behind this paper is then the following formula,
which is a special case of a far more general result known as the
Gaussian kinematic
formula. (cf.\ \cite{RFG,Taylor2006,Taylor:Adler:2003} for details.)
It says that for Gaussian related $f=F(g)$ field defined over well-behaving sets
$S$,
\beq
\label{ECApproximation}
\E \left\{\p \left(\left\{t\in S\: f(t) \geq u\right\}\right) \right\}
\ =\ \sum_{j=0}^{\dim S} \lips_{j}(S) \, \rho_j^F(u).  
\eeq 
The $\lips_j(S)$ are the \LK\ curvatures, or intrinsic volumes, of $S$,
and estimating these in statistical practice was one of the aims of
this paper. If the $g^i$ are stationary and isotropic, and the
variance of their first order partial derivatives is given by the second
spectral moment $\lambda_2$, then $\lips_j(S) = \lambda_2^{j/2}
\lips^E_j(S)$ where the $\lips^E_j(S)$ are the Euclidean \LKCs.  The
Euclidean \LKCs are known under a variety of names, including
Quermassintegrales, Minkowski or Steiner functionals, integral
curvatures, and intrinsic volumes, the differences between them
generally being of ordering and scaling. In general, $\lips^E_j(S)$
can be thought of as a measure of the `$j$-dimensional size' of $S$.
For example, when $N = 2$, $\lips^E_2(S)$ is the two dimensional area of
$S$, $\lips^E_1(S)$ is half its boundary length, and $\lips^E_0(S)$
its Euler characteristic.  When $N=3$, $\lips^E_3(S)$ is the
three-dimensional volume of $S$, $\lips^E_2(S)$ is half the surface
area, $\lips^E_1(S)$ is twice the caliper diameter of $S$, (where the
caliper diameter of a convex $S$ is defined by placing the solid
between two parallel planes (or calipers), measuring the distance
between the planes, and averaging over all rotations of $S$) and
$\lips^E_0(S)$ is again the \EC.  In fact $\lips^E_0(S)=\lips_0(S)=\p
(S)$ in all cases.

While these are simple examples of \LK\ curvatures, when the underlying $g^i$
are not isotropic their covariance induces a Riemannian metric on $S$,
and the $\lips_j(S)$ become quite complicated. To be a little more 
specific, they
involve integrals of traces of powers of the curvature tensor and
second fundamental forms, with respect to the volume form determined
by the induced metric. In brief, these are not typically things that a
practitioner wants to compute analytically, particularly since
computing them also requires knowledge of the covariance function of
the $g^i$, something which typically needs to be estimated from data.

On the other hand, the functions $\rho^F_j$ in
(\ref{ECApproximation}) are much easier to compute, and in many cases are given by
\beqq
\rho_j^F(u) \ =\  
(-1)^j (2\pi)^{-j/2}\frac{d^j}{dx^j}P\left\{F(Z)\geq x \right\}
\Big|_{x=g(u)},
\eeqq
where $Z\sim N(0,I_{k\times k})$ and $g$ is a function determined by $F$.
In such cases, computing the $\rho_k$ is thus simply a matter
of calculus, and for most interesting $F$ this has already been done.
(\citet*{RFG,ARF} have a number of useful examples with references to
others.) A particularly simple example is given in the Gaussian case,
in which, as mentioned above, $k=1$ and $F$ is the identity
function. Then the $\rho_j$ are then given by  
\beqq
\rho_j(u)\ =\
 { (2\pi)^{-(j+1)/2} } H_{j-1}(u)e^{-u^2/2} \qquad j=0,1,\dots N,
\eeqq
where $H_j$, $j\geq 1$ is the $j$-th Hermite polynomial and 
\beqq
H_{-1}(u)
\definedas \sqrt{2\pi}e^{u^2/2}P(N(0,1)\geq u).
\eeqq

To give an example in a more complicated case, we consider the $\chi^2$ case,
with $k$ degrees of freedom, so that $F$ is given by the leftmost  function  in 
\eqref{rob-chi}. In this case, the $\rho_j$ are, for
$j \geq 1$ and $u>0$, 
\beqq
\lefteqn{\rho_{j}(u) = \frac{u^{(k-j)/2}e^{-u/2}}{(2 \pi)^{j/2} \Gamma(k/2) 2^{(k-2)/2}}\sum_{l=0}^{\lfloor \frac{j-1}{2} \rfloor} \sum_{m=0}^{j-1-2l} } \\
&& \qquad \quad\times   1_{\{k \geq j-m-2l\}} \binom{k-1}{j-1-m-2l} \frac{(-1)^{j-1+m+l}(j-1)!}{m! l! 2^l}  u^{m+l}.
\eeqq
When $j=0$, 
\beqq
\rho_{0}(u) \ = \  \P\left\{\chi^2_k \geq u\right\}.
\eeqq
While this may look a little complicated, it is trivial to code, and the LKC regression procedure then continues as in the Gaussian case.

\bibliographystyle{imsart-nameyear}
\bibliography{paper.bib}

\begin{thebibliography}{40}

\bibitem[\protect\citeauthoryear{Adler}{1981}]{GRF}
\begin{bbook}[author]
\bauthor{\bsnm{Adler},~\bfnm{R.~J.}\binits{R.~J.}}
(\byear{1981}).
\btitle{{The Geometry of Random Fields}}.
\bpublisher{John Wiley \& Sons Ltd.}, \baddress{Chichester}.
\bnote{Reprinted in 2010 by SIAM}.
\bmrnumber{82h:60103}
\end{bbook}
\endbibitem

\bibitem[\protect\citeauthoryear{Adler}{2000}]{Adler:2000}
\begin{barticle}[author]
\bauthor{\bsnm{Adler},~\bfnm{R.~J.}\binits{R.~J.}}
(\byear{2000}).
\btitle{{On Excursion Sets, Tube Formulae, and Maxima of Random Fields}}.
\bjournal{The Annals of Applied Probability}
\bvolume{10}
\bpages{1--74}.
\end{barticle}
\endbibitem

\bibitem[\protect\citeauthoryear{Adler, Subag and Taylor}{2012}]{AST}
\begin{barticle}[author]
\bauthor{\bsnm{Adler},~\bfnm{R.~J.}\binits{R.~J.}},
  \bauthor{\bsnm{Subag},~\bfnm{E.}\binits{E.}} \AND
  \bauthor{\bsnm{Taylor},~\bfnm{J.~E.}\binits{J.~E.}}
(\byear{2012}).
\btitle{{Rotation and Scale Space Random Fields and the Gaussian Kinematic
  Formula}}.
\bjournal{Annals of Statistics}
\bvolume{40}
\bpages{2910--2942}.
\end{barticle}
\endbibitem

\bibitem[\protect\citeauthoryear{Adler and Taylor}{2007}]{RFG}
\begin{bbook}[author]
\bauthor{\bsnm{Adler},~\bfnm{R.~J.}\binits{R.~J.}} \AND
  \bauthor{\bsnm{Taylor},~\bfnm{J.~E.}\binits{J.~E.}}
(\byear{2007}).
\btitle{Random Fields and Geometry}.
\bpublisher{Springer}.
\end{bbook}
\endbibitem

\bibitem[\protect\citeauthoryear{Adler and Taylor}{2011}]{ATSF}
\begin{bbook}[author]
\bauthor{\bsnm{Adler},~\bfnm{R.~J.}\binits{R.~J.}} \AND
  \bauthor{\bsnm{Taylor},~\bfnm{J.~E.}\binits{J.~E.}}
(\byear{2011}).
\btitle{Topological Complexity of Smooth Random Functions}.
\bseries{Lecture Notes in Mathematics}
\bvolume{2019}.
\bpublisher{Springer}, \baddress{Heidelberg}.
\bnote{Lectures from the 39th Probability Summer School held in Saint-Flour,
  2009, {\'E}cole d'{\'E}t{\'e} de Probabilit{\'e}s de Saint-Flour.
  [Saint-Flour Probability Summer School]}.
\bmrnumber{2768175}
\end{bbook}
\endbibitem

\bibitem[\protect\citeauthoryear{Adler, Taylor and Worsley}{2015?}]{ARF}
\begin{bbook}[author]
\bauthor{\bsnm{Adler},~\bfnm{R.~J.}\binits{R.~J.}},
  \bauthor{\bsnm{Taylor},~\bfnm{J.~E.}\binits{J.~E.}} \AND
  \bauthor{\bsnm{Worsley},~\bfnm{K.~J.}\binits{K.~J.}}
(\byear{2015}?).
\btitle{Applications of Random Fields and Geometry: Foundations and Case
  Studies}.
\bpublisher{Springer-Verlag}
\bnote{In preparation, early chapters available at
  http://webee.technion.ac.il/people/adler/publications.html}.
\end{bbook}
\endbibitem

\bibitem[\protect\citeauthoryear{Bardeen et~al.}{1986}]{Bardeen}
\begin{barticle}[author]
\bauthor{\bsnm{Bardeen},~\bfnm{J.~M.}\binits{J.~M.}},
  \bauthor{\bsnm{Bond},~\bfnm{J.~R.}\binits{J.~R.}},
  \bauthor{\bsnm{Kaiser},~\bfnm{N.}\binits{N.}} \AND
  \bauthor{\bsnm{Szalay},~\bfnm{A.~S.}\binits{A.~S.}}
(\byear{1986}).
\btitle{{The Statistics of Peaks of Gaussian Random Fields}}.
\bjournal{Astrophysical Journal}
\bvolume{304}
\bpages{15--61}.
\end{barticle}
\endbibitem

\bibitem[\protect\citeauthoryear{Beckmann, Jenkinson and
  Smith}{{2003}}]{Beckmann2003}
\begin{barticle}[author]
\bauthor{\bsnm{Beckmann},~\bfnm{C.~F.}\binits{C.~F.}},
  \bauthor{\bsnm{Jenkinson},~\bfnm{M.}\binits{M.}} \AND
  \bauthor{\bsnm{Smith},~\bfnm{S.~M.}\binits{S.~M.}}
(\byear{{2003}}).
\btitle{{General Multilevel Linear Modeling for Group Analysis in FMRI}}.
\bjournal{{Neuroimage}}
\bvolume{{20}}
\bpages{{1052-1063}}.
\end{barticle}
\endbibitem

\bibitem[\protect\citeauthoryear{Cao and
  Worsley}{1999}]{Cao:Worsley:1999:Correlation}
\begin{barticle}[author]
\bauthor{\bsnm{Cao},~\bfnm{J.}\binits{J.}} \AND
  \bauthor{\bsnm{Worsley},~\bfnm{K.~J.}\binits{K.~J.}}
(\byear{1999}).
\btitle{{The Geometry of Correlation Fields With an Application to Functional
  Connectivity of the Brain}}.
\bjournal{Annals of Applied Probability}
\bvolume{9}
\bpages{1021--1057}.
\bmrnumber{2000k:60099}
\end{barticle}
\endbibitem

\bibitem[\protect\citeauthoryear{Dehaene-Lambertz et~al.}{{2006}}]{Dehaene2006}
\begin{barticle}[author]
\bauthor{\bsnm{Dehaene-Lambertz},~\bfnm{G}\binits{G.}},
  \bauthor{\bsnm{Dehaene},~\bfnm{S}\binits{S.}},
  \bauthor{\bsnm{Anton},~\bfnm{J.~L}\binits{J.~L.}},
  \bauthor{\bsnm{Campagne},~\bfnm{A}\binits{A.}},
  \bauthor{\bsnm{Ciuciu},~\bfnm{P}\binits{P.}},
  \bauthor{\bsnm{Dehaene},~\bfnm{G.~P.}\binits{G.~P.}},
  \bauthor{\bsnm{Denghien},~\bfnm{I}\binits{I.}},
  \bauthor{\bsnm{Jobert},~\bfnm{A}\binits{A.}},
  \bauthor{\bsnm{LeBihan},~\bfnm{D}\binits{D.}},
  \bauthor{\bsnm{Sigman},~\bfnm{M}\binits{M.}},
  \bauthor{\bsnm{Pallier},~\bfnm{C}\binits{C.}} \AND
  \bauthor{\bsnm{Poline},~\bfnm{J.~B}\binits{J.~B.}}
(\byear{{2006}}).
\btitle{{Functional Segregation of Cortical Language Areas by Sentence
  Repetition}}.
\bjournal{{Human Brain Mapping}}
\bvolume{{27}}
\bpages{{360-371}}.
\bdoi{{10.1002/hbm.20250}}
\end{barticle}
\endbibitem

\bibitem[\protect\citeauthoryear{Fantaye et~al.}{2014}]{Marinucci-LKC}
\begin{barticle}[author]
\bauthor{\bsnm{Fantaye},~\bfnm{Yabebal}\binits{Y.}},
  \bauthor{\bsnm{Hansen},~\bfnm{Frode}\binits{F.}},
  \bauthor{\bsnm{Maino},~\bfnm{Davide}\binits{D.}} \AND
  \bauthor{\bsnm{Marinucci},~\bfnm{Domenico}\binits{D.}}
(\byear{2014}).
\btitle{Cosmological Applications of the Gaussian Kinematic Formula}.
\bnote{arxiv.org/pdf/1406.5420v1.pdf}.
\end{barticle}
\endbibitem

\bibitem[\protect\citeauthoryear{Friston et~al.}{1994}]{FRISTON94}
\begin{barticle}[author]
\bauthor{\bsnm{Friston},~\bfnm{K.~J.}\binits{K.~J.}},
  \bauthor{\bsnm{Worsley},~\bfnm{K.~J.}\binits{K.~J.}},
  \bauthor{\bsnm{Frackowiak},~\bfnm{R.~S.~J.}\binits{R.~S.~J.}},
  \bauthor{\bsnm{Mazziotta},~\bfnm{J.~C.}\binits{J.~C.}} \AND
  \bauthor{\bsnm{Evans},~\bfnm{A.~C.}\binits{A.~C.}}
(\byear{1994}).
\btitle{{Assessing the Significance of Focal Activations Using Their Spatial
  Extent}}.
\bjournal{Human Brain Mapping}
\bvolume{1}
\bpages{214--220}.
\end{barticle}
\endbibitem

\bibitem[\protect\citeauthoryear{Gott et~al.}{2007}]{gott2007genus}
\begin{barticle}[author]
\bauthor{\bsnm{Gott},~\bfnm{J.~R.}\binits{J.~R.}},
  \bauthor{\bsnm{Colley},~\bfnm{W.~N.}\binits{W.~N.}},
  \bauthor{\bsnm{Park},~\bfnm{C.~G.}\binits{C.~G.}},
  \bauthor{\bsnm{Park},~\bfnm{C}\binits{C.}} \AND
  \bauthor{\bsnm{Mugnolo},~\bfnm{C}\binits{C.}}
(\byear{2007}).
\btitle{{Genus Topology of the Cosmic Microwave Background From the WMAP 3-year
  Data}}.
\bjournal{Monthly Notices of the Royal Astronomical Society}
\bvolume{377}
\bpages{1668--1678}.
\end{barticle}
\endbibitem

\bibitem[\protect\citeauthoryear{Gott et~al.}{2008}]{Sloan2008}
\begin{barticle}[author]
\bauthor{\bsnm{Gott},~\bfnm{J.~R.}\binits{J.~R.}},
  \bauthor{\bsnm{Hambrick},~\bfnm{D.~C.}\binits{D.~C.}},
  \bauthor{\bsnm{Vogeley},~\bfnm{M.~S.}\binits{M.~S.}},
  \bauthor{\bsnm{Kim},~\bfnm{J}\binits{J.}},
  \bauthor{\bsnm{Park},~\bfnm{C}\binits{C.}},
  \bauthor{\bsnm{Choi},~\bfnm{Y.~Y.}\binits{Y.~Y.}},
  \bauthor{\bsnm{Cen},~\bfnm{R}\binits{R.}} \AND
  \bauthor{\bsnm{Ostriker},~\bfnm{K}\binits{K.} \bsuffix{J.~P.~Nagamine}}
(\byear{2008}).
\btitle{Genus Topology of Structure in the Sloan Digital Sky Survey: Model
  Testing}.
\bjournal{Astrophysical Journal}
\bvolume{675}
\bpages{16--28}.
\end{barticle}
\endbibitem

\bibitem[\protect\citeauthoryear{Hall, Fisher and Hoffman}{{1994}}]{Hall1994}
\begin{barticle}[author]
\bauthor{\bsnm{Hall},~\bfnm{P}\binits{P.}},
  \bauthor{\bsnm{Fisher},~\bfnm{N.~I.}\binits{N.~I.}} \AND
  \bauthor{\bsnm{Hoffman},~\bfnm{B}\binits{B.}}
(\byear{{1994}}).
\btitle{{On the Nonparametric Estimation of Covariance Functions}}.
\bjournal{{Annals of Statistics}}
\bvolume{{22}}
\bpages{{2115-2134}}.
\end{barticle}
\endbibitem

\bibitem[\protect\citeauthoryear{Imiya and Eckhardt}{{1999}}]{Imiya1999}
\begin{barticle}[author]
\bauthor{\bsnm{Imiya},~\bfnm{A}\binits{A.}} \AND
  \bauthor{\bsnm{Eckhardt},~\bfnm{U}\binits{U.}}
(\byear{{1999}}).
\btitle{{The Euler Characteristics of Discrete Objects and Discrete
  Quasi-objects}}.
\bjournal{{Computer Vision and Image Understanding}}
\bvolume{{75}}
\bpages{{307-318}}.
\end{barticle}
\endbibitem

\bibitem[\protect\citeauthoryear{Kilner and Friston}{2010}]{Friston-Kilner}
\begin{barticle}[author]
\bauthor{\bsnm{Kilner},~\bfnm{J.~M.}\binits{J.~M.}} \AND
  \bauthor{\bsnm{Friston},~\bfnm{K.~J.}\binits{K.~J.}}
(\byear{2010}).
\btitle{{Topological Inference for EEG and MEG Data}}.
\bjournal{Annals of Applied Statistics}
\bvolume{4}
\bpages{1272-1290}.
\end{barticle}
\endbibitem

\bibitem[\protect\citeauthoryear{Loader}{1999}]{Loader1999}
\begin{bbook}[author]
\bauthor{\bsnm{Loader},~\bfnm{C}\binits{C.}}
(\byear{1999}).
\btitle{Local Regression and Likelihood}.
\bpublisher{Springer}, \baddress{New York}.
\end{bbook}
\endbibitem

\bibitem[\protect\citeauthoryear{Loader}{2010}]{Loader2010}
\begin{bmanual}[author]
\bauthor{\bsnm{Loader},~\bfnm{C}\binits{C.}}
(\byear{2010}).
\btitle{locfit: Local Regression, Likelihood and Density Estimation}
\bnote{R package version 1.5-6}.
\end{bmanual}
\endbibitem

\bibitem[\protect\citeauthoryear{Mantoglou and Wilson}{{1982}}]{Mantoglou1982}
\begin{barticle}[author]
\bauthor{\bsnm{Mantoglou},~\bfnm{A}\binits{A.}} \AND
  \bauthor{\bsnm{Wilson},~\bfnm{J.~L.}\binits{J.~L.}}
(\byear{{1982}}).
\btitle{{The Turning Bands Method for Simulation of Random Fields Using Line
  Generation by a Spectral Method}}.
\bjournal{{Water Resources Research}}
\bvolume{{18}}
\bpages{{1379-1394}}.
\end{barticle}
\endbibitem

\bibitem[\protect\citeauthoryear{Pukelsheim}{2006}]{Pukelsheim2006}
\begin{bbook}[author]
\bauthor{\bsnm{Pukelsheim},~\bfnm{F}\binits{F.}}
(\byear{2006}).
\btitle{Optimal Design of Experiments}.
\bpublisher{Society for Industrial and Applied Mathematics}.
\end{bbook}
\endbibitem

\bibitem[\protect\citeauthoryear{Rao}{{1962}}]{Rao1962}
\begin{barticle}[author]
\bauthor{\bsnm{Rao},~\bfnm{C.~R.}\binits{C.~R.}}
(\byear{{1962}}).
\btitle{{A Note on a Generalized Inverse of a Matrix With Applications to
  Problems in Mathematical-statistics}}.
\bjournal{{Journal of the Royal Statistical Society Series B-Statistical
  Methodology}}
\bvolume{{24}}
\bpages{{152-158}}.
\end{barticle}
\endbibitem

\bibitem[\protect\citeauthoryear{Sampson and Guttorp}{{1992}}]{Sampson1992}
\begin{barticle}[author]
\bauthor{\bsnm{Sampson},~\bfnm{PD}\binits{P.}} \AND
  \bauthor{\bsnm{Guttorp},~\bfnm{P}\binits{P.}}
(\byear{{1992}}).
\btitle{{Nonparametric Estimation of Nonstationary Spatial Covariance
  Structure}}.
\bjournal{{Journal of the American Statistical Association}}
\bvolume{{87}}
\bpages{{108-119}}.
\end{barticle}
\endbibitem

\bibitem[\protect\citeauthoryear{Schlather}{2009}]{Schlather2009}
\begin{bmanual}[author]
\bauthor{\bsnm{Schlather},~\bfnm{M}\binits{M.}}
(\byear{2009}).
\btitle{RandomFields: Simulation and Analysis of Random Fields}
\bnote{R package version 1.3.41}.
\end{bmanual}
\endbibitem

\bibitem[\protect\citeauthoryear{Shafie
  et~al.}{2003}]{Shafie:Sigal:Siegmund:Worsley:2001}
\begin{barticle}[author]
\bauthor{\bsnm{Shafie},~\bfnm{K.}\binits{K.}},
  \bauthor{\bsnm{Sigal},~\bfnm{B.}\binits{B.}},
  \bauthor{\bsnm{Siegmund},~\bfnm{D.}\binits{D.}} \AND
  \bauthor{\bsnm{Worsley},~\bfnm{K.~J.}\binits{K.~J.}}
(\byear{2003}).
\btitle{{Rotation Space Random Fields With an Application to fMRI Data}}.
\bjournal{Annals of Statistics}
\bvolume{31}
\bpages{1732--1771}.
\end{barticle}
\endbibitem

\bibitem[\protect\citeauthoryear{Siegmund and Worsley}{1995}]{S-WOR95}
\begin{barticle}[author]
\bauthor{\bsnm{Siegmund},~\bfnm{D.~O.}\binits{D.~O.}} \AND
  \bauthor{\bsnm{Worsley},~\bfnm{K.~J.}\binits{K.~J.}}
(\byear{1995}).
\btitle{Testing for a signal with unknown location and scale in a stationary
  {G}aussian random field}.
\bjournal{Ann. Statist}
\bvolume{23}
\bpages{608--639}.
\end{barticle}
\endbibitem

\bibitem[\protect\citeauthoryear{Smith et~al.}{{2004}}]{Smith2004}
\begin{barticle}[author]
\bauthor{\bsnm{Smith},~\bfnm{S.~M.}\binits{S.~M.}},
  \bauthor{\bsnm{Jenkinson},~\bfnm{M}\binits{M.}},
  \bauthor{\bsnm{Woolrich},~\bfnm{M.~W.}\binits{M.~W.}},
  \bauthor{\bsnm{Beckmann},~\bfnm{C.~F.}\binits{C.~F.}},
  \bauthor{\bsnm{Behrens},~\bfnm{T.~E.~J.}\binits{T.~E.~J.}},
  \bauthor{\bsnm{Johansen-Berg},~\bfnm{H}\binits{H.}},
  \bauthor{\bsnm{Bannister},~\bfnm{P.~R.}\binits{P.~R.}},
  \bauthor{\bsnm{De~Luca},~\bfnm{M}\binits{M.}},
  \bauthor{\bsnm{Drobnjak},~\bfnm{I}\binits{I.}},
  \bauthor{\bsnm{Flitney},~\bfnm{D.~E.}\binits{D.~E.}},
  \bauthor{\bsnm{Niazy},~\bfnm{R.~K.}\binits{R.~K.}},
  \bauthor{\bsnm{Saunders},~\bfnm{J}\binits{J.}},
  \bauthor{\bsnm{Vickers},~\bfnm{J}\binits{J.}},
  \bauthor{\bsnm{Zhang},~\bfnm{Y.~Y.}\binits{Y.~Y.}},
  \bauthor{\bsnm{De~Stefano},~\bfnm{N}\binits{N.}},
  \bauthor{\bsnm{Brady},~\bfnm{J.~M.}\binits{J.~M.}} \AND
  \bauthor{\bsnm{Matthews},~\bfnm{P.~M.}\binits{P.~M.}}
(\byear{{2004}}).
\btitle{{Advances in Functional and Structural MR Image Analysis and
  Implementation as FSL}}.
\bjournal{{Neuroimage}}
\bvolume{{23}}
\bpages{{S208-S219}}.
\end{barticle}
\endbibitem

\bibitem[\protect\citeauthoryear{Taylor}{{2006}}]{Taylor2006}
\begin{barticle}[author]
\bauthor{\bsnm{Taylor},~\bfnm{J.~E.}\binits{J.~E.}}
(\byear{{2006}}).
\btitle{{A Gaussian kinematic formula}}.
\bjournal{{Annals of Probability}}
\bvolume{{34}}
\bpages{{122-158}}.
\end{barticle}
\endbibitem

\bibitem[\protect\citeauthoryear{Taylor and Adler}{2003}]{Taylor:Adler:2003}
\begin{barticle}[author]
\bauthor{\bsnm{Taylor},~\bfnm{J.~E.}\binits{J.~E.}} \AND
  \bauthor{\bsnm{Adler},~\bfnm{R.~J.}\binits{R.~J.}}
(\byear{2003}).
\btitle{Euler characteristics for {G}aussian fields on manifolds}.
\bjournal{Annals of Probability}
\bvolume{31}
\bpages{533--563}.
\end{barticle}
\endbibitem

\bibitem[\protect\citeauthoryear{Taylor, Takemura and
  Adler}{2005}]{Taylor:Takemura:Adler:2003:Validity}
\begin{barticle}[author]
\bauthor{\bsnm{Taylor},~\bfnm{J.~E.}\binits{J.~E.}},
  \bauthor{\bsnm{Takemura},~\bfnm{A}\binits{A.}} \AND
  \bauthor{\bsnm{Adler},~\bfnm{R.~J.}\binits{R.~J.}}
(\byear{2005}).
\btitle{{Validity of the Expected Euler Characteristic Heuristic}}.
\bjournal{Annals of Probability}
\bvolume{33}
\bpages{1362--1396}.
\bmrnumber{MR2150192 (2006b:60073)}
\end{barticle}
\endbibitem

\bibitem[\protect\citeauthoryear{Taylor and
  Worsley}{{2006}}]{Taylor2006Inference}
\begin{barticle}[author]
\bauthor{\bsnm{Taylor},~\bfnm{J.~E.}\binits{J.~E.}} \AND
  \bauthor{\bsnm{Worsley},~\bfnm{K.~J.}\binits{K.~J.}}
(\byear{{2006}}).
\btitle{{Inference for Magnitudes and Delays of Responses in the FIAC Data
  Using BRAINSTAT/FMRISTAT}}.
\bjournal{{Human Brain Mapping}}
\bvolume{{27}}
\bpages{{434-441}}.
\bdoi{{10.1002/hbm.20248}}
\end{barticle}
\endbibitem

\bibitem[\protect\citeauthoryear{Taylor and Worsley}{{2007}}]{Taylor2007}
\begin{barticle}[author]
\bauthor{\bsnm{Taylor},~\bfnm{J.~E.}\binits{J.~E.}} \AND
  \bauthor{\bsnm{Worsley},~\bfnm{K.~J.}\binits{K.~J.}}
(\byear{{2007}}).
\btitle{{Detecting Sparse Signals in Random Fields, With an Application to
  Brain Mapping}}.
\bjournal{{Journal of the American Statistical Association}}
\bvolume{{102}}
\bpages{{913-928}}.
\bdoi{{10.1198/016214507000000815}}
\end{barticle}
\endbibitem

\bibitem[\protect\citeauthoryear{Theis and Tanaka}{{2005}}]{Theis2005}
\begin{bincollection}[author]
\bauthor{\bsnm{Theis},~\bfnm{F.~J.}\binits{F.~J.}} \AND
  \bauthor{\bsnm{Tanaka},~\bfnm{T}\binits{T.}}
(\byear{{2005}}).
\btitle{{A fast and efficient method for compressing fMRI data sets}}.
In \bbooktitle{{Artificial Neural Networks: Formal Models and Their
  Applications - ICANN 2005, Pt. 2, Proceedings}},
(\beditor{~\bsnm{{Duch, W and Kacprzyk, J and Oja, E and Zadrozny, S}}}, ed.).
\bseries{{Lecture Notes in Computer Science}}
\bvolume{{3697}}
\bpages{{769-777}}.
\end{bincollection}
\endbibitem

\bibitem[\protect\citeauthoryear{Torres}{1994}]{TORRES}
\begin{barticle}[author]
\bauthor{\bsnm{Torres},~\bfnm{S}\binits{S.}}
(\byear{1994}).
\btitle{{Topological Analysis of COBE-DMR Cosmic Microwave Background Maps}}.
\bjournal{Astrophysical Journal}
\bvolume{423}
\bpages{9--12}.
\end{barticle}
\endbibitem

\bibitem[\protect\citeauthoryear{Vogeley et~al.}{1994}]{VPGHG}
\begin{barticle}[author]
\bauthor{\bsnm{Vogeley},~\bfnm{M.~S.}\binits{M.~S.}},
  \bauthor{\bsnm{Park},~\bfnm{C}\binits{C.}},
  \bauthor{\bsnm{Geller},~\bfnm{M.~J.}\binits{M.~J.}},
  \bauthor{\bsnm{Huchin},~\bfnm{J.~P.}\binits{J.~P.}} \AND
  \bauthor{\bsnm{Gott},~\bfnm{J.~R.}\binits{J.~R.}}
(\byear{1994}).
\btitle{{Topological Analysis of the CfA Redshift Survey}}.
\bjournal{Astrophysical Journal}
\bvolume{420}
\bpages{525--544}.
\end{barticle}
\endbibitem

\bibitem[\protect\citeauthoryear{Wilde et~al.}{{2009}}]{Wilde2009}
\begin{barticle}[author]
\bauthor{\bsnm{Wilde},~\bfnm{Michael}\binits{M.}},
  \bauthor{\bsnm{Foster},~\bfnm{Ian}\binits{I.}},
  \bauthor{\bsnm{Iskra},~\bfnm{Kamil}\binits{K.}},
  \bauthor{\bsnm{Beckman},~\bfnm{Pete}\binits{P.}},
  \bauthor{\bsnm{Zhang},~\bfnm{Zhao}\binits{Z.}},
  \bauthor{\bsnm{Espinosa},~\bfnm{Allan}\binits{A.}},
  \bauthor{\bsnm{Hategan},~\bfnm{Mihael}\binits{M.}},
  \bauthor{\bsnm{Clifford},~\bfnm{Ben}\binits{B.}} \AND
  \bauthor{\bsnm{Raicu},~\bfnm{Ioan}\binits{I.}}
(\byear{{2009}}).
\btitle{{Parallel Scripting for Applications at the Petascale and Beyond}}.
\bjournal{{Computer}}
\bvolume{{42}}
\bpages{{50-60}}.
\end{barticle}
\endbibitem

\bibitem[\protect\citeauthoryear{Worsley}{1994}]{Worsley:1994:Chi:t:F}
\begin{barticle}[author]
\bauthor{\bsnm{Worsley},~\bfnm{K.~J.}\binits{K.~J.}}
(\byear{1994}).
\btitle{{Local Maxima and the Expected Euler Characteristic of Excursion Sets
  of $\chi\sp 2,\ {F}$ and $t$ Fields}}.
\bjournal{Advances in Applied Probability}
\bvolume{26}
\bpages{13--42}.
\bmrnumber{94i:60064}
\end{barticle}
\endbibitem

\bibitem[\protect\citeauthoryear{Worsley}{1995a}]{Worsley:1995:Boundary}
\begin{barticle}[author]
\bauthor{\bsnm{Worsley},~\bfnm{K.~J.}\binits{K.~J.}}
(\byear{1995}a).
\btitle{{Boundary Corrections for the Expected Euler Characteristic of
  Excursion Sets of Random Fields, With an Application to Astrophysics}}.
\bjournal{Advances in Applied Probability}
\bvolume{27}
\bpages{943--959}.
\bmrnumber{97b:60068}
\end{barticle}
\endbibitem

\bibitem[\protect\citeauthoryear{Worsley}{1995b}]{WOR95a}
\begin{barticle}[author]
\bauthor{\bsnm{Worsley},~\bfnm{K.~J.}\binits{K.~J.}}
(\byear{1995}b).
\btitle{{Estimating the Number of Peaks in a Random Field Using the Hadwiger
  Characteristic of Excursion Sets, With Applications to Medical Images}}.
\bjournal{Annals of Statistics}
\bvolume{23}
\bpages{640--669}.
\end{barticle}
\endbibitem

\bibitem[\protect\citeauthoryear{Zhao et~al.}{2007}]{Zhao2007}
\begin{barticle}[author]
\bauthor{\bsnm{Zhao},~\bfnm{Yong}\binits{Y.}},
  \bauthor{\bsnm{Hategan},~\bfnm{Mihael}\binits{M.}},
  \bauthor{\bsnm{Clifford},~\bfnm{Ben}\binits{B.}},
  \bauthor{\bsnm{Foster},~\bfnm{Ian}\binits{I.}}, \bauthor{\bparticle{von}
  \bsnm{Laszewski},~\bfnm{Gregor}\binits{G.}},
  \bauthor{\bsnm{Nefedova},~\bfnm{Veronika}\binits{V.}},
  \bauthor{\bsnm{Raicu},~\bfnm{Ioan}\binits{I.}},
  \bauthor{\bsnm{Stef-Praun},~\bfnm{Tiberiu}\binits{T.}} \AND
  \bauthor{\bsnm{Wilde},~\bfnm{Michael}\binits{M.}}
(\byear{2007}).
\btitle{{Swift: Fast, Reliable, Loosely Coupled Parallel Computation}}.
\bjournal{{IEEE Congress on Services}}
\bvolume{0}
\bpages{199-206}.
\end{barticle}
\endbibitem

\end{thebibliography}

\address{Robert J.\ Adler \\Electrical Engineering\\
Technion, Haifa,
Israel 32000 \\
\printead{e1}\\
\printead{u1}}

\noindent  \address{Kevin Bartz\\
Renaissance Technologies LLC\\
600 Route 25A, 
East Setauket, NY 11733\\
\printead{e2}
}

\noindent  \address{S.\ C.\ Kou\\
Department of Statistics\\ Harvard University\\
1 Oxford Street     
Cambridge, MA 02138 \\
\printead{e3}\\
\printead{u3}}

\address{Anthea Monod\\
Department of Systems Biology\\ Columbia University\\
1130 St. Nicholas Avenue\\
New York, NY 10032\\
\printead{e4}\\
\printead{u4}}

\end{document}